\documentclass[12pt,a4paper]{article}
\usepackage{a4wide}
\usepackage{amsfonts,amsmath}
\usepackage{latexsym}
\usepackage[nobysame,initials]{amsrefs}
\allowdisplaybreaks
\usepackage{times}
\title{Strengthened volume inequalities for $L_p$ zonoids of even isotropic measures }
\author{K\'aroly J. B\"or\"oczky, Ferenc Fodor, Daniel Hug}

\newcommand{\proof}{\noindent{\it Proof: }}
\newcommand{\proofbox}{\mbox{ $\Box$}\\}
\newcommand{\R}{\mathbb{R}}
\newcommand{\Z}{\mathbb{Z}}

\newcommand{\NN}{\text{\sf N}}

\newtheorem{lemma}{Lemma}[section]
\newtheorem{theo}[lemma]{Theorem}
\newtheorem{claim}[lemma]{Claim}

\newtheorem{coro}[lemma]{Corollary}
\newtheorem{conj}[lemma]{Conjecture}
\newtheorem{prop}[lemma]{Proposition}

\begin{document}

\maketitle

\begin{abstract}
We strengthen the volume inequalities for $L_p$ zonoids of even isotropic measures and 
for their duals, which are due to Ball, Barthe and Lutwak, Yang, Zhang. Along the way, 
we prove a stronger version of the Brascamp-Lieb inequality for a family of functions 
that can approximate arbitrary well some Gaussians when equality holds. 
The special case $p=\infty$ yields a stability version of the reverse isoperimetric 
inequality for centrally symmetric bodies.
\end{abstract}

\section{Introduction}

According to the classical isoperimetric inequality Euclidean balls minimize the 
surface area among convex bodies of given volume in Euclidean space $\R^n$. We call 
a subset of $\R^n$ a convex body if it is compact, convex and has non-empty interior. 
Let $B^n$ be the Euclidean unit ball centred at the origin, and let $S(\cdot)$ and 
$V(\cdot)$ denote the surface area and the volume functional in $\R^n$, respectively. 
The isoperimetric inequality can be stated in the form
$$
\frac{S(B^n)^n}{V(B^n)^{n-1}}\le \frac{S(K)^n}{V(K)^{n-1}},
$$
where equality holds if and only if $K$ is a Euclidean ball. Recently, N. Fusco, 
F. Maggi, A. Pratelli \cite{FMP08}  proved  an essentially optimal stability version 
of the isoperimetric inequality. It states that if $K$ is a convex body with $V(K)=V(B^n)$ and 
if $S(B^n)\geq (1-\varepsilon)S(K)$ holds
for some small $\varepsilon>0$, then $K$ is close to some
translate  $B^n+x$, $x\in\R^n$, of the unit ball; namely,
$$
V(K\Delta(B^n+x))\leq \gamma\varepsilon^{1/2},
$$
where $\gamma>0$ depends only on $n$, and $\Delta$ denotes the symmetric difference of sets.

Stability estimates for the planar isoperimetric inequality
go back to the works of Minkowski and Bonnesen. However, 
a systematic exploration  is much more recent.
We refer to the survey articles of H.~Groemer  \cite{Groemer1990, Groemer1993} for an 
introduction to geometric stability results. The recent monograph \cite{Sch14} by 
R.~Schneider provides an up-to-date treatment of the topic including applications.
Here we only note that the stability estimate related to the isoperimetric inequality 
obtained in \cite{FMP08} 
was extended to a stability version of the Brunn-Minkowski inequality by
 A.~Figalli, F.~Maggi, A.~Pratelli \cite{FMP09, FMP10}.

Aiming at a reverse isoperimetric inequality, F.~Behrend \cite{Behrend1937} suggested 
to consider equivalence classes of convex bodies with respect to non-singular linear 
transformations. C.M.~Petty \cite{Petty1961} proved
(see also A. Giannopoulos, M. Papadimitrakis \cite{GianPapa1999}) that there is an 
essentially unique  representative minimizing the isoperimetric ratio in each equivalence class. 
The unique minimizer in an equivalence class is characterized by the property that
its suitably normalized area measure is isotropic. We give a precise definition of 
isotropic measures later.  
This characterization yields that cubes minimize the isoperimetric
ratio within the class of parallelotopes, and regular simplices within the class of simplices.

The functional that assigns to each equivalence class the minimum of the isoperimetric
ratio within that class is affine invariant and upper semi-continuous, therefore it attains its 
maximum on the affine equivalence classes of convex bodies.
In the Euclidean plane, the method of
F.~Behrend \cite{Behrend1937} 
yields that the maximum is attained by the affine equivalence class of triangles, and by the 
affine equivalence class of parallelograms if the convex body is assumed to be centrally symmetric.
The extension of these results to higher dimensions proved to be quite difficult. 
Decades after Behrend's paper, K.M.~Ball in \cite{Bal89, Bal91a} managed to establish reverse 
forms of the isoperimetric inequality in arbitrary dimensions. 
More precisely, the largest isoperimetric ratio is attained
by simplices according to \cite{Bal91a}, and by parallelotopes among centrally symmetric 
convex bodies according to \cite{Bal89}. Since the reverse isoperimetric inequality and a 
stronger form of it for general convex bodies are discussed  
in K.J. B\"or\"oczky, D. Hug \cite{BoH}, in this paper we concentrate on centrally symmetric 
convex bodies.

In order to state the result of K.M.~Ball \cite{Bal89} about centrally symmetric convex bodies, 
we set $W^n=[-1,1]^n$, and note that
$S(W^n)=n2^n=nV(W^n)$.

\bigskip

\noindent{\bf Theorem~A } (K.M.~Ball) {\it For any centrally symmetric convex body $K$ 
in $\R^n$, there exists some
$\Phi\in{\rm GL}(n)$ such that}
\begin{equation}\label{riiec}
\frac{S(\Phi K)^n}{V(\Phi K)^{n-1}}\leq \frac{S(W^n)^n}{V(W^n)^{n-1}}.
\end{equation}

\bigskip

The case of equality in Theorem~A was settled by F.~Barthe \cite{Bar98}. 
He proved that if the left side of 
\eqref{riiec} 
is minimized over all $\Phi\in{\rm GL}(n)$, then equality holds precisely when 
$K$ is a  parallelotope.

Our first objective is to prove a stability version of the reverse isoperimetric 
inequality for centrally symmetric convex bodies. Following
 \cite{FMP09, FMP10, FMP08}, we define an affine invariant distance
of origin symmetric convex bodies $K$ and $M$ based on the volume difference. Let $\alpha=V(K)^{-{1}/n}$,
$\beta=V(M)^{-{1}/n}$, and define
$$
\delta_{\rm vol}(K,M)=\min\left\{V\left(\Phi (\alpha K)\Delta (\beta M)\right):\,
\Phi\in{\rm SL}(n)\right\}
$$
where ${\rm SL}(n)$ is the group of linear transformations of $\R^n$ of determinant one.
In fact, $\delta_{\rm vol}(\cdot,\cdot)$ induces a metric on the linear equivalence classes 
of origin symmetric convex bodies.

The John ellipsoid of a convex body $K$ in $\R^n$ is the unique maximum volume ellipsoid 
contained in $K$. If $K$ is origin symmetric, then its
John ellipsoid is also origin symmetric. Note that each convex body
has an affine image whose John ellipsoid is $B^n$. The John ellipsoid is
a frequently used  tool in geometric analysis, and, in particular, 
it was used by K.M.~Ball in the proof of the reverse isoperimetric inequality. Since we 
will use the John ellipsoid in our arguments, below we review its basic properties 
(see (\ref{John10})). For a more detailed treatment of the topic, we refer to K.M.~Ball \cite{Bal03}, 
P.M.~Gruber \cite{Gru07} and R.~Schneider \cite{Sch14}.

\begin{theo}
\label{vol}
Let $K$ be an origin symmetric convex body in $\R^n$, $n\ge 3$, whose John ellipsoid is 
a Euclidean ball, and let $\varepsilon\in [0,1)$. 
If $\delta_{\rm vol}(K,W^n)\geq \varepsilon$, then
$$
\frac{S(K)^n}{V(K)^{n-1}}\leq (1-\gamma\,\varepsilon^3)\frac{S(W^n)^n}{V(W^n)^{n-1}},
$$
where  $\gamma=n^{-cn^3}$ for some absolute constant $c>0$.
\end{theo}

The stability order (the exponent $3$ of $\varepsilon$) in Theorem~\ref{vol} is
close to be optimal, but most probably it is not optimal. Considering a convex body 
$K$ which is obtained from $W^n$ by cutting off  simplices of height $\varepsilon$ 
at  the vertices of $W^n$, one can see that the exponent of $\varepsilon$ must be at least $1$ in Theorem~\ref{vol}.

Another common affine invariant distance between
convex bodies is the  Banach-Mazur metric $\delta_{\rm BM}(K,M)$, which we 
define here only for origin symmetric convex bodies $K$ and $M$. Let
$$
\delta_{\rm BM}(K,M)=\log\min\{\lambda\geq 1:\,
K\subseteq \Phi(M)\subseteq \lambda \,K\mbox{ for some }
\Phi\in{\rm GL}(n)\}.
$$
We note that $\delta_{\rm vol}\le 2n^2\delta_{\rm BM}$ (see, say, \cite{BoH}). 
Furthermore,
$\delta_{\rm BM}\le \gamma \ \delta_{\rm vol}^{\frac{1}{n}}$, where 
$\gamma$ depends only on the dimension $n$ 
(see \cite[Section 5]{BHenk}). The example of a ball from which a cap 
is cut off shows that in the latter inequality the exponent $\frac1n$ 
cannot be replaced by anything larger than $\frac2{n+1}$.

\begin{theo}
\label{BM}
Let $K$ be an origin symmetric convex body in $\R^n$, $n\ge 3$, whose 
John ellipsoid is a Euclidean ball, and let $\varepsilon\in [0,1)$. 
If $\delta_{\rm BM}(K,W^n)\geq \varepsilon$, then
$$
\frac{S(K)^n}{V(K)^{n-1}}\leq (1-\gamma\,\varepsilon^n)\frac{S(W^n)^n}{V(W^n)^{n-1}},
$$
where $\gamma=n^{-cn^3}$ for some absolute constant $c>0$.
\end{theo}

The stability order (the exponent $n$ of $\varepsilon$) in Theorem~\ref{BM} is again
close to be optimal, but very likely it is not optimal. Considering a convex 
body $K$ which is obtained from $W^n$ by cutting off  simplices of height 
$\varepsilon$ at  the vertices of $W^n$, one can see that the exponent of 
$\varepsilon$ must be at least $n-1$ in Theorem~\ref{BM}. 

In the planar case, a modification of the argument of F.~Behrend 
\cite{Behrend1937} leads to stability results of optimal order.

\begin{theo}
\label{volBMplanar}
Let $K$ be an origin symmetric convex body in $\R^2$ which has a square as an inscribed parallelogram of maximum area. 
Let  $\varepsilon\in [0,1)$. 
If $\delta_{\rm vol}(K,W^2)\geq \varepsilon$ or 
$\delta_{\rm BM}(K,W^2)\geq \varepsilon$, then
$$
\frac{S(K)^2}{V(K)}\leq \left(1-\frac{\varepsilon}{54}\right)\frac{S(W^2)^2}{V(W^2)}.
$$
\end{theo}

Note that for an origin symmetric convex body $K$ in $\R^2$ there always exists a linear transform $\Phi\in \text{GL}(2)$ 
such that a square is an inscribed parallelogram of maximum area of $\Phi K$. 
In particular, if we define $\text{ir}(K)=\min\{S(\Phi K)^2/V(\Phi K):\Phi\in \text{GL}(2)\}$,  
for an origin symmetric convex body in $K$ in $\R^2$, and if $\varepsilon\in [0,1)$, then Theorem \ref{volBMplanar} 
implies that
$$
\text{ir}(K)\leq \left(1-\frac{\varepsilon}{54}\right)\text{ir}(W^2)
$$
provided that $\delta_{\rm vol}(K,W^2)\geq \varepsilon$ or $\delta_{\rm BM}(K,W^2)\geq \varepsilon$.
\medskip

As mentioned before, the proof of the reverse isoperimetric inequality 
by K.M.~Ball \cite{Bal89, Bal91a} is based on a volume estimate for convex bodies
whose John ellipsoid is the unit ball $B^n$. Let $S^{n-1}$ denote the Euclidean unit sphere.
 According to a classical theorem of F.~John \cite{Joh37} (see also K.M.~Ball \cite{Bal03}),
$B^n$ is the ellipsoid of maximal volume in an origin symmetric convex body 
$K$ if and only if $B^n\subseteq K$  and there exist
$\pm u_1,\ldots,\pm u_k\in S^{n-1}\cap \partial K$ and $c_1,\ldots,c_k>0$ such that
\begin{equation}
\label{John10}
\sum_{i=1}^kc_i u_i\otimes u_i={\rm Id}_n,
\end{equation}
where $\otimes$ denotes the tensor product of vectors in $\R^n$, ${\rm Id}_n$ 
denotes the $n\times n$ identity matrix and $\partial K$ is the boundary of $K$.

Following A.~Giannopoulos, M.~Papadimitrakis \cite{GianPapa1999} and
E.~Lutwak, D.~Yang, G.~Zhang \cite{Lutwak1}, we call
an even Borel measure $\mu$ on the unit sphere $S^{n-1}$ isotropic if
$$
\int_{S^{n-1}}u\otimes u\, d\mu(u)={\rm Id}_n.
$$
In this case, equating traces of both sides we obtain that $\mu(S^{n-1})=n$.

Using the standard notation $\langle\cdot\,,\cdot\rangle$ for the Euclidean scalar product  and 
$\|\cdot\|$ for the induced norm in $\R^n$,  
the support function $h_K$ of a convex compact set $K$ in $\R^n$ at $v\in \R^n$ is defined as 
$$
h_K(v)=\max\{\langle v,x\rangle:\,x\in K\}.
$$
For any $p\ge 1$ and an even measure $\mu$ on $S^{n-1}$ not concentrated on any great subsphere, 
we define the $L_p$ zonoid $Z_p(\mu)$ associated with $\mu$ by
$$
h_{Z_p(\mu)}(v)^p=\int_{S^{n-1}}|\langle u,v\rangle|^p\,d\mu(u),
$$
which is  a zonoid in the classical sense if $p=1$.
In addition, let
$$
Z_\infty(\mu)=\lim_{p\to\infty}Z_p(\mu)={\rm conv}\,{\rm supp}\,\mu,
$$
and for $1\leq p\leq\infty$, let $Z_p^*(\mu)$ be the polar of $Z_p(\mu)$.
In particular,
\begin{eqnarray*}
Z_p^*(\mu)&=&\left\{x\in\R^n:\,\int_{S^{n-1}}|\langle x,u\rangle|^p\,d\mu(u)\leq 1
\right\} \mbox{ \ for $p\in[1,\infty)$},\\
Z_\infty^*(\mu)&=&\{x\in\R^n:\,\langle x,u\rangle\leq 1 \mbox{ for }u\in{\rm supp}\,\mu\},
\end{eqnarray*}
and hence $Z_2(\mu)=B^n$ for any even isotropic measure $\mu$.

It follows from D.R. Lewis \cite{Lew78} (see also
E. Lutwak, D. Yang and G. Zhang \cite{Lutwak0,LYZ05}) that any $n$-dimensional 
subspace of $L_p$ is isometric to
$\|\cdot\|_{Z_p^*(\mu)}$ for some isotropic measure $\mu$ on $S^{n-1}$, where
$$
\|x\|_{Z_p^*(\mu)}=\left(\int_{S^{n-1}}|\langle x,u\rangle|^p\,d\mu(u)\right)^{\frac1p},\qquad x\in\R^n.
$$

We call a measure $\nu$ on $S^{n-1}$  a cross measure if there is an orthonormal 
basis $u_1,\ldots,u_n$ of $\R^n$ such that
$$
{\rm supp}\,\nu=\{\pm u_1,\ldots,\pm u_n\},
$$
and $\nu(\{u_i\})=\nu(\{-u_i\})=1/2$ for $i=1,\ldots,n$, and hence $\nu$ is even and isotropic. 
We fix a cross measure $\nu_n$ on $S^{n-1}$.
We note that if
$p\in[1,\infty]$, and $\Gamma(\cdot)$ is Euler's Gamma function, then
$$
V(Z_p(\nu_n))=
\left\{\begin{array}{ll}
\frac{\Gamma(1+\frac{n}2)\Gamma(1+\frac{p}2)}
{\Gamma(1+\frac{1}2)\Gamma(1+\frac{n+p}2)}
&\mbox{ \ if }p\geq 1,\\[
2ex]
\frac{2^n}{n!}&\mbox{ \ if }p=\infty.
\end{array}
\right.
$$
In addition,
$$
V(Z^*_p(\nu_n))=
\left\{\begin{array}{ll}
2^n\frac{\Gamma(1+\frac1{p})^n}{\Gamma(1+\frac{n}{p})}
&\mbox{ \ if }p\geq 1,
\\[2ex]
2^n&\mbox{ \ if }p=\infty.
\end{array}
\right.
$$
The crucial statement leading to the reverse isoperimetric inequality is the case of $Z^*_\infty(\mu)$.

\bigskip

\noindent{\bf Theorem~B } {\it If $\mu$ is an even isotropic measure on $S^{n-1}$ and
$p\in[1,\infty]$, then}
\begin{eqnarray*}
\label{Zp}
V(Z_p(\mu))&\geq&V(Z_p(\nu_n)), \\
\label{Zp*}
V(Z_p^*(\mu))&\leq& V(Z^*_p(\nu_n)).
\end{eqnarray*}
{\it Assuming $p\neq 2$, equality holds if and only if $\mu$ is a cross measure.}

\bigskip

Theorem~B is the work of K.M. Ball \cite{Bal91a} and F. Barthe \cite{Bar98} 
if $\mu$ is discrete, and their method  was extended to arbitrary even isotropic 
measures $\mu$ by E. Lutwak, D. Yang, and G. Zhang \cite{Lutwak0}. The measures 
on $S^{n-1}$ which have an isotropic linear image are characterized by
K.J.~B\"or\"oczky, E.~Lutwak, D.~Yang and  G.~Zhang \cite{BLYZ15}, building on the works of
E.A.~Carlen,
and D.~Cordero-Erausquin \cite{CCE09},  J.~Bennett, A.~Carbery, M.~Christ and 
T.~Tao \cite{BCCT08} and B.~Klartag \cite{Kla10}.
We note that isotropic measures on $\R^n$ play a central role in the
KLS conjecture by R.~Kannan, L.~Lov\'asz and   M.~Simonovits
\cite{KLM95}; see, for instance,
F.~Barthe and  D.~Cordero-Erausquin \cite{BCE13}, O.~Guedon and  E.~Milman \cite{GuM11} and
B.~Klartag \cite{Kla09}.

To state a stability version of Theorem~B, a natural notion of distance between 
two isotropic measures $\mu$ and $\nu$ is the Wasserstein distance (also called 
the Kantorovich-Monge-Rubinstein distance) $\delta_W(\mu,\nu)$. To define it, we write
$\angle(v,w)$ to denote the angle between non-zero vectors $v$ and $w$; that is, 
the geodesic distance of the unit vectors
 $\|v\|^{-1}v$ and $\|w\|^{-1}w$ on the unit sphere. Let
 ${\rm Lip}_1(S^{n-1})$ denote the family of Lipschitz functions with Lipschitz 
constant at most $1$; 
namely, $f:\,S^{n-1}\to\R$ is in ${\rm Lip}_1(S^{n-1})$ if 
$\|f(x)-f(y)\|\leq \angle(x,y)$ for
$x,y\in S^{n-1}$. Then the Wasserstein distance of $\mu$ and $\nu$ is given by
$$
\delta_W(\mu,\nu)=\max\left\{\int_{S^{n-1}}f\,d\mu-\int_{S^{n-1}}f\,d\nu:\,f\in 
{\rm Lip}_1(S^{n-1})\right\}.
$$
What we actually need in this paper is the Wasserstein distance of an isotropic 
measure $\mu$ from the closest cross measure. Therefore,
in the case of two isotropic measures $\mu$ and $\nu$, we define
$$
\delta_{\rm WO}(\mu,\nu)=\min\left\{\delta_W(\mu,\Phi_*\nu):\,\Phi\in {\rm O}(n)\right\}
$$
where $\Phi_*\nu$ denotes the pushforward of $\nu$ by $\Phi: S^{n-1}\to S^{n-1}$.

\begin{theo}
\label{Zpmustab}
Let $\mu$ be an even isotropic measure on $S^{n-1}$, $n\ge 2$,  let $\varepsilon\in[0,1)$,
and let $p\in[1,\infty]$ with $p\neq 2$. If $\delta_{\rm WO}(\mu,\nu_n)\geq\varepsilon>0$, then
\begin{eqnarray*}
V(Z_p(\mu))&\ge& (1+\gamma\varepsilon^3) V(Z_p(\nu_n)),\\
V(Z^*_p(\mu))&\le& (1-\gamma\varepsilon^3) V(Z^*_p(\nu_n))
\end{eqnarray*}
where $\gamma=n^{-cn^3}\min\{|p-2|^2,1\}$ for an absolute constant $c>0$.
\end{theo}

To state another stability version of Theorem~B, in the case  $p=\infty$,  
we use the ``spherical" Hausdorff distance 
$\delta_H(X,Y)$ of compact  sets $X,Y\subseteq S^{n-1}$ given by
$$
\delta_H(X,Y)=\min\left\{\max_{x\in X}\min_{y\in Y}\angle(x,y),\max_{y\in Y}\min_{x\in X}\angle(x,y)\right\}.
$$
In addition, let
$$
\delta_{HO}(X,Y)=\min\left\{\delta_H(X,\Phi Y):\,\Phi\in {\rm O}(n)\right\}.
$$

We note that if $\delta_{HO}({\rm supp}\,\mu,{\rm supp}\,\nu_n)\leq 1/(7n^2)$
for an even isotropic measure $\mu$, then $\delta_{WO}(\mu,\nu_n)\leq 2n 
\delta_{HO}({\rm supp}\,\mu,{\rm supp}\,\nu_n)$  according to 
Corollary~\ref{Hausdorffcross}. However, as we will see in 
Section~\ref{secreverseisostab}, Theorem~\ref{Zpmustab} 
implies the following seemingly stronger statement in the case $p=\infty$.

\begin{coro}
\label{Zinftymustab}
If  $\mu$ is an even isotropic measure on $S^{n-1}$, and 
$\delta_{HO}({\rm supp}\,\mu,{\rm supp}\,\nu_n)\geq \varepsilon>0$, then
\begin{eqnarray*}
V(Z_{\infty}(\mu))&\ge& (1+\gamma\varepsilon^3) V(Z_{\infty}(\nu_n)),\\
V(Z^*_{\infty}(\mu))&\le& (1-\gamma\varepsilon^3) V(Z^*_{\infty}(\nu_n))
\end{eqnarray*}
where $\gamma=n^{-cn^3}$ for an absolute constant $c>0$.
\end{coro}


We note that the order $\varepsilon^3$ of the error term in Corollary~\ref{Zinftymustab} can be improved to 
$\varepsilon$ if $n=2$ according to Theorem~\ref{Zinftymustab2}.


The proof of Theorem~B by is based on the rank one case of the geometric Brascamp-Lieb inequality. 
An essential tool in  our approach is the proof provided
by  F.~Barthe \cite{Bar97,Bar98}, which is based on mass transportation. 
Therefore, we review the argument from \cite{Bar97} in Section~\ref{secbrascamp-lieb}. 
At the end of that section,  we outline the arguments leading to
Theorem~\ref{vol}, Theorem \ref{BM} and
Theorem~\ref{Zpmustab}  and we describe the structure of the paper. We also indicate 
in Section~\ref{secbrascamp-lieb} what stability result can be expected concerning the Brascamp-Lieb 
inequality (see Conjecture~\ref{BL-RBLsamefunction}). Along the way of proving our main statements, 
we also establish some properties of arbitrary (not only even) isotropic measures in Section~\ref{seciso}
 that might be useful in other applications as well.

Let us point out that the corresponding question in the non-symmetric setting is wide open.
We call an isotropic measure $\mu$ on $S^{n-1}$ centred if
$$\int_{S^{n-1}}u\, d\mu (u)=o.$$
Here and in the following, we write $o$ for the origin (the zero vector). 
 For a centred isotropic measure $\mu$ on $S^{n-1}$, and for $p\in[1,\infty)$,
we define the non-symmetric $L_p$ zonoid  $Z_p(\mu)$ by
\begin{eqnarray*}
h_{Z_p(\mu)}(v)^p&=&2\int_{S^{n-1}}\max\{0,\langle v,u\rangle\}^p\,d\mu(u),\\
Z_p^*(\mu)&=&\left\{x\in\R^n:\,\int_{S^{n-1}}\max\{0,\langle x,u\rangle\}^p\,d\mu(u)\leq \frac12
\right\}.
\end{eqnarray*}
This notion (for any discrete measure on $S^{n-1}$, not only isotropic ones), occurs in 
M. Weberndorfer \cite{Web13} in connection with reverse versions of the Blaschke-Santal\'o inequality.
The factor $2$ is included to match the earlier definition for even isotropic measures.
The difference to the case of even isotropic measures is that if $p=2$ and $\mu$ is a non-even 
centered isotropic measure, then  $Z_2(\mu)$ is typically not a Euclidean ball but has 
constant squared width; namely,
$h_{Z_p(\mu)}(v)^2+h_{Z_p(\mu)}(-v)^2$ is constant for $v\in S^{n-1}$.

\begin{conj}
\label{conj:nonsym}
If  $\mu$ is a centered isotropic measure  on $S^{n-1}$ and  $p\in[1,\infty)$, moreover
$\nu$ is an isotropic measure on $S^{n-1}$ such that 
${\rm supp}\,\nu$ consists of the vertices of a regular simplex, then
\begin{eqnarray*}
\label{Zpconj}
V(Z_p(\mu))&\geq &V(Z_p(\nu)),\\
\label{Zp*conj}
V(Z_p^*(\mu))&\leq & V(Z_p^*(\nu)).
\end{eqnarray*}
\end{conj}

If  $\mu$ is a centered isotropic measure  on $S^{n-1}$, then 
$Z_\infty(\mu)={\rm conv\;supp}\,\mu$. In particular, if $p=\infty$,
then (\ref{Zp*conj}) was proved by K.M. Ball in \cite{Bal91a} for discrete $\mu$, 
 (\ref{Zpconj})  was proved by F. Barthe in \cite{Bar98} again for discrete $\mu$, 
and the case of general centered isotropic $\mu$ was handled 
E. Lutwak, D. Yang and G. Zhang \cite{Lutwak1}.

An inequality related to the case $p=2$ of Conjecture~\ref{conj:nonsym} is proved by
E.~Lutwak, D.~Yang, G.~Zhang \cite{LYZ10}.

\section{A brief review of the Brascamp-Lieb and the reverse Brascamp-Lieb  inequality}
\label{secbrascamp-lieb}

The rank one geometric Brascamp-Lieb inequality (\ref{BL}), identified by  K.M. Ball \cite{Bal89}
as an essential case
of the rank one Brascamp-Lieb inequality,
due to H.J.~Brascamp, E.H.~Lieb \cite{BrL76}, and the reverse form (\ref{RBL}), due to F. Barthe \cite{Bar97,Bar98},
read as follows.
If $u_1,\ldots,u_k\in S^{n-1}$ are distinct unit vectors and $c_1,\ldots,c_k>0$ satisfy
$$
\sum_{i=1}^kc_i u_i\otimes u_i={\rm Id}_n,
$$
and $f_1,\ldots,f_k$ are non-negative measurable functions on $\R$, then
\begin{eqnarray}
\label{BL}
\int_{\R^n}\prod_{i=1}^kf_i(\langle x,u_i\rangle)^{c_i}\,dx&\leq&
\prod_{i=1}^k\left(\int_{\R}f_i\right)^{c_i}, \mbox{ \ and}\\
\label{RBL}
\int_{\R^n}^{\ast}\sup_{x=\sum_{i=1}^kc_i\theta_iu_i}\prod_{i=1}^kf_i(\theta_i)^{c_i}\,dx&\geq&
\prod_{i=1}^k\left(\int_{\R}f_i\right)^{c_i}.
\end{eqnarray}
In (\ref{RBL}), the supremum extends over all $\theta_1,\ldots,\theta_k\in\R$. Since the integrand 
need not be a measurable function, 
we have to consider the outer integral. If $k=n$, then $u_1,\ldots,u_n$ form an orthonormal basis 
and therefore $\theta_1,\ldots,\theta_k$ are uniquely determined for a given $x\in\R^n$.

According to  F.~Barthe \cite{Bar98}, if equality holds in (\ref{BL}) or in (\ref{RBL}) and none 
of the functions $f_i$ is identically zero or a scaled version of a  Gaussian, then there is an 
origin symmetric regular crosspolytope in $\R^n$ such that $u_1,\ldots,u_k$ lie among its vertices. 
Conversely, equality holds in (\ref{BL}) and (\ref{RBL}) if  each $f_i$ is a scaled 
version of the same centered Gaussian,
or if $k=n$ and $u_1,\ldots,u_n$ form an  orthonormal basis.

A thorough discussion of the rank one Brascamp-Lieb inequality can be found in
E.~Carlen, D.~Cordero-Erausquin \cite{CCE09}. The higher rank case,
due to E.H.~Lieb \cite{Lie90}, is reproved and further explored by F.~Barthe \cite{Bar98} 
(including a discussion of the equality case), and is again carefully analysed by
J.~Bennett, T.~Carbery, M.~Christ, T.~Tao \cite{BCCT08}. In particular, see F.~Barthe, 
D.~Cordero-Erausquin, M.~Ledoux, 
B.~Maurey \cite{BCLM11} for an enlightening review of the relevant literature and an 
approach via Markov semigroups in a quite general framework.

 F.~Barthe \cite{Bar97, Bar98}  provided concise proofs of (\ref{BL}) and (\ref{RBL})
based on mass transportation (see also K.M.~Ball \cite{Bal03} for  (\ref{BL})). We sketch 
the main ideas of his  approach, since it will be the starting point of subsequent refinements.

We assume that
each $f_i$ is a positive continuous probability density both for (\ref{BL}) and (\ref{RBL}), 
and let $g(t)=e^{-\pi t^2}$ be the Gaussian density.
For $i=1,\ldots,k$, we consider the transportation map
$T_i:\R\to\R$ satisfying
$$
\int_{-\infty}^tf_i(s)\,ds=\int_{-\infty}^{T_i(t)} g(s)\,ds.
$$
It is easy to see that $T_i$ is bijective, differentiable and
\begin{equation}
\label{masstrans}
f_i(t)=g(T_i(t))\cdot T'_i(t),\qquad t\in\R.
\end{equation}
To these transportation maps, we associate the smooth transformation $\Theta:\R^n\to\R^n$ given by
$$
\Theta(x)=\sum_{i=1}^kc_iT_i(\langle u_i,x\rangle )\,u_i,\qquad x\in\R^n,
$$
which satisfies
$$
d\Theta(x)=\sum_{i=1}^kc_iT'_i(\langle u_i,x\rangle )\,u_i\otimes u_i.
$$
In this case, $d\Theta(x)$ is positive definite and $\Theta:\R^n\to\R^n$ is
 injective (see \cite{Bar97,Bar98}). We will need the following two estimates due to K.M.~Ball \cite{Bal89}
(see also \cite{Bar98} for a simpler proof of (i)).
\begin{description}
\item{(i)} For any $t_1,\ldots,t_k>0$, we have
$$
\det \left(\sum_{i=1}^kt_ic_i u_i\otimes u_i\right)\geq \prod_{i=1}^k t_i^{c_i}.
$$

\item{(ii)}
If $z=\sum_{i=1}^kc_i\theta_i u_i$ for
$\theta_1,\ldots,\theta_k\in\R$, then
\begin{equation}
\label{RBLfunc}
\|z\|^2\le \sum_{i=1}^kc_i\theta_i^2.
\end{equation}
\end{description}
Therefore, using first \eqref{masstrans},  then (i) with $t_i=T'_i(\langle u_i,x\rangle)$, 
 the definition of $\Theta$ and (ii), and finally the transformation formula, 
 the following argument leads to the Brascamp-Lieb inequality (\ref{BL}).
\begin{align}
\label{BLstep1}
\int_{\R^n}\prod_{i=1}^kf_i(\langle u_i,x\rangle)^{c_i}\,dx&=
\int_{\R^n}\left(\prod_{i=1}^kg(T_i(\langle u_i,x\rangle))^{c_i}\right)
\left(\prod_{i=1}^kT'_i(\langle u_i,x\rangle)^{c_i}\right)\,dx\\
\label{BLstep2}
&\leq \int_{\R^n}\left(\prod_{i=1}^ke^{-\pi c_iT_i(\langle u_i,x\rangle)^2}\right)
 \det\left(\sum_{i=1}^kc_iT'_i(\langle u_i,x\rangle )\,u_i\otimes u_i\right)\,dx\\
\nonumber
&\leq  \int_{\R^n}e^{-\pi \|\Theta(x)\|^2}\det\left( d\Theta(x)\right)\,dx\\
\nonumber
&\leq  \int_{\R^n}e^{-\pi \|y\|^2}\,dy=1.
\end{align}
The Brascamp-Lieb inequality (\ref{BL}) for arbitrary non-negative integrable 
functions $f_i$ follows by scaling and approximation.

For the reverse Brascamp-Lieb inequality (\ref{RBL}), we consider the inverse $S_i$ of $T_i$, and hence
\begin{eqnarray}
\nonumber
\int_{-\infty}^tg(s)\,ds=\int_{-\infty}^{S_i(t)} f_i(s)\,ds,\\
\label{masstransS}
g(t)=f_i(S_i(t))\cdot S'_i(t),\qquad t\in\R.
\end{eqnarray}
In addition,
$$
d\Psi(x)=\sum_{i=1}^kc_iS'_i(\langle u_i,x\rangle )\,u_i\otimes u_i
$$
holds for the smooth transformation $\Psi:\R^n\to\R^n$ given by
$$
\Psi(x)=\sum_{i=1}^kc_iS_i(\langle u_i,x\rangle )\,u_i,\qquad x\in\R^n.
$$
In particular, $d\Psi(x)$ is positive definite and $\Psi:\R^n\to\R^n$ is
 injective (see \cite{Bar97,Bar98}).
Therefore (i) and (\ref{masstransS}) lead to
\begin{align}
\nonumber
&\int_{\R^n}^{\ast}\sup_{x=\sum_{i=1}^kc_i\theta_iu_i}\prod_{i=1}^kf_i(\theta_i)^{c_i}\,dx\nonumber\\
&\qquad \geq
\int_{\R^n}^{\ast}\left(\sup_{\Psi(y)=\sum_{i=1}^kc_i\theta_iu_i}\prod_{i=1}^kf_i(\theta_i)^{c_i}\right)
\det\left( d\Psi(y)\right)\,dy\nonumber\\
\label{RBLstep1}
&\qquad\geq  \int_{\R^n}\left(\prod_{i=1}^kf_i(S_i(\langle u_i,y\rangle))^{c_i} \right)
\det\left(\sum_{i=1}^kc_iS'_i(\langle u_i,y\rangle )\,u_i\otimes u_i\right)\,dy\\
\label{RBLstep2}
&\qquad\geq \int_{\R^n}\left(\prod_{i=1}^kf_i(S_i(\langle u_i,y\rangle))^{c_i}\right)
\left(\prod_{i=1}^kS'_i(\langle u_i,y\rangle)^{c_i}\right)\,dy\\
\nonumber
&\qquad = \int_{\R^n}\left(\prod_{i=1}^kg(\langle u_i,y\rangle)^{c_i}\right)
\,dy=  \int_{\R^n}e^{-\pi \|y\|^2}\,dy=1.
\end{align}
Again, the reverse Brascamp-Lieb inequality (\ref{RBL}) for arbitrary non-negative integrable functions $f_i$ follows by scaling and approximation.

We observe that (i) shows that the optimal constant in the geometric Brascamp-Lieb inequality is $1$.
The stability version  of (i) (with $v_i=\sqrt{c_i} u_i$), Lemma~\ref{Ball-Barthe-stab}, is an essential tool in proving
a stability version of the Brascamp-Lieb inequality leading to Theorem~\ref{Zpmustab}.

Even if we do not use it in this paper, we point out that F. Barthe \cite{Bar04} proved ``continuous" versions of the Brascamp-Lieb and the reverse Brascamp-Lieb inequalities that work for any isotropic measure $\mu$ on $S^{n-1}$
(see (\ref{contBL}) and (\ref{contRBL}) below). Here we only consider the case in which all non-negative real functions involved coincide with a ``nice" probability density function, which is the common case in geometric applications. So let 
$f:\R\to[0,\infty)$ be such that $\int_R f=1$ and 
${\rm supp}(f)=[a,b]$   for some $a,b\in[-\infty,\infty]$. Further, we assume that $f$ is positive and continuous on $[a,b]$.  According to 
\cite{Bar04}, we have
\begin{equation}
\label{contBL}
\int_{\R^n}\exp\left(\int_{S^{n-1}}\log f(\langle x,u\rangle)\,d\mu(u)\right)\,dx\leq 1.
\end{equation}
For the reverse inequality, let $h:\R^n\to[0,\infty)$ be a measurable function which satisfies  
$$
h\left(\int_{S^{n-1}}\theta(u)\, u\,d\mu(u)\right)\geq
\exp\left(\int_{S^{n-1}}\log f(\theta(u))\,d\mu(u)\right)
$$
for any continuous function $\theta:{\rm supp}\,\mu\to \R$. Then, we have
\begin{equation}
\label{contRBL}
\int_{\R^n}h\geq 1.
\end{equation}

Let us briefly discuss how K.M.~Ball \cite{Bal89} and  F. Barthe \cite{Bar98} used the Brascamp-Lieb inequality and its reverse form to prove
 the discrete version of Theorem~B.  In this section, we write $\mu$ to denote the isotropic measure on $S^{n-1}$ whose support is $\{u_1,\ldots,u_k\}$ with $\mu(\{u_i\})=c_i$, and we assume that $\mu$ is an even measure. 
For $i=1,\ldots,k$, we consider the probability densities on $\R$ (see (\ref{tpint})) given by
$$
f_i(t)=\frac1{2\Gamma(1+\frac{1}p)}\,e^{-|t|^p},\qquad t\in\R,
$$
if $p\in[1,\infty)$, and $f_i=\frac12\mathbf{1}_{[-1,1]}$ if $p=\infty$, where
$$
\mathbf{1}_{[-1,1]}(t)=\left\{
\begin{array}{rl}
1 &\mbox{ if }t\in[-1,1],\\
0&\mbox{ otherwise}.
\end{array}\right.
$$
We will frequently use the following observation due to K. Ball \cite{Bal91a}. If $K$ is an orgin symmetric convex body in $\R^n$ with associated norm $\|\cdot\|_K$ and if $p\in[1,\infty)$, then
$$
V(K)=\frac1{\Gamma(1+\frac{n}p)}\int_{\R^n}e^{-\|x\|_K^p}\,dx,
$$
where 
$$
\|x\|_K=\min\{\lambda\geq 0:\,x\in\lambda K\},\qquad x\in\R^n.
$$
In particular, if $p\in[1,\infty)$,  then
\begin{eqnarray}
\nonumber
V(Z_p^*(\mu))&=&\frac1{\Gamma(1+\frac{n}p)}\int_{\R^n}\exp\left(-\sum_{i=1}^kc_i|\langle x,u_i\rangle|^p\right)\,dx\\
\label{volZ*p0}
&=&
\frac{2^n\Gamma\left(1+\frac1p \right)^n}{\Gamma(1+\frac{n}p)}\int_{\R^n}\prod_{i=1}^kf_i(\langle x,u_i\rangle)^{c_i}\,dx\\
&\leq& \frac{2^n\Gamma\left(1+\frac1p \right)^n}{\Gamma(1+\frac{n}p)}
\prod_{i=1}^k\left(\int_{\R}f_i\right)^{c_i}=\frac{2^n\Gamma\left(1+\frac1p \right)^n}{\Gamma(1+\frac{n}p)}.
\label{volZ*p}
\end{eqnarray}
On the other hand, if $p=\infty$, then using $f_i=\tfrac{1}{2}\mathbf{1}_{[-1,1]}$, we have
$$
V(Z_\infty^*(\mu))=2^n\int_{\R^n}\prod_{i=1}^kf_i(\langle x,u_i\rangle)^{c_i}\,dx\leq
2^n\prod_{i=1}^k\left(\int_{\R}f_i\right)^{c_i}=2^n.
$$
Equality in (\ref{volZ*p}) leads to equality in the Brascamp-Lieb inequality, and hence $k=2n$ and
$u_1,\ldots,u_k$ form the vertices of a regular crosspolytope in $\R^n$.

For the lower bound on the volume of the $L_p$ zonotopes and $p\in[1,\infty]$, let us choose $p^*\in[1,\infty]$ 
such that $\frac1p+\frac1{p^*}=1$. If $p\in[1,\infty)$, then an (auxiliary) origin symmetric convex body is defined by
$$
M_p(\mu)=\left\{\sum_{i=1}^kc_i\theta_iu_i:\, \sum_{i=1}^kc_i|\theta_i|^p\leq 1\right\}.
$$
We drop the reference to $\mu$, if it does not cause any misunderstanding. In particular,
$$
\|x\|_{M_p}=\left(\inf_{x=\sum_{i=1}^kc_i\theta_iu_i}  \sum_{i=1}^kc_i|\theta_i|^p\right)^{\frac1p},\qquad x\in\R^n.
$$
In addition, we define
$$
M_\infty(\mu)=\left\{\sum_{i=1}^kc_i\theta_iu_i:\, |\theta_i|\leq 1\mbox{ for }i=1,\ldots,k\right\}.
$$

We claim that if $p\in[1,\infty]$, then
\begin{equation}
\label{MpZp*}
M_p(\mu)\subseteq Z_{p^*}(\mu).
\end{equation}
Let $x\in M_p(\mu)$, and hence $x=\sum_{i=1}^kc_i\theta_iu_i$ with 
$\sum_{i=1}^kc_i|\theta_i|^p\leq 1$ if $p\in[1,\infty)$,
and $|\theta_i|\leq 1$ for $i=1,\ldots,k$ if $p=\infty$. If $p\in(1,\infty)$, 
then it follows from H\"older's inequality that, for any $v\in\R^n$, we have
$$
\langle x,v\rangle=\sum_{i=1}^kc_i\theta_i\langle u_i,v\rangle\leq
\left(\sum_{i=1}^kc_i|\theta_i|^p\right)^{\frac1p}\left(\sum_{i=1}^kc_i|\langle u_i,v\rangle|^{p^*}\right)^{\frac1{p^*}}
\leq h_{Z_{p^*}}(v).
$$
If $p=1$, then
$$
\langle x,v\rangle=\sum_{i=1}^kc_i\theta_i\langle u_i,v\rangle\leq \max_{i=1,\ldots,k}|\langle u_i,v\rangle|
=h_{Z_{\infty}}(v).
$$
In addition, if $p=\infty$, then
$$
\langle x,v\rangle=\sum_{i=1}^kc_i\theta_i\langle u_i,v\rangle\leq \sum_{i=1}^kc_i|\langle u_i,v\rangle|=
h_{Z_1}(v).
$$

Now if  $p\in[1,\infty)$, then we deduce from (\ref{MpZp*}) and the reverse Brascamp-Lieb inequality (\ref{RBL}) that
\begin{eqnarray}
\nonumber
V(Z_{p^*}(\mu))&\geq &V(M_p(\mu))=\frac1{\Gamma(1+\frac{n}p)}\int_{\R^n}\exp\left(-\|x\|_{M_p}^p\right)\,dx\\
\label{volZp*RBL}
&=&\frac{2^n\Gamma(1+\frac{1}p)^n}{\Gamma(1+\frac{n}p)}\int_{\R^n}^{\ast}\sup_{x=\sum_{i=1}^kc_i\theta_iu_i}  \prod_{i=1}^kf_i(\theta_i)^{c_i}\,dx\\
\label{volZp*}
&\geq&\frac{2^n\Gamma(1+\frac{1}p)^n}{\Gamma(1+\frac{n}p)}\prod_{i=1}^k\left(\int_{\R}f_i\right)^{c_i}=
\frac{2^n\Gamma(1+\frac{1}p)^n}{\Gamma(1+\frac{n}p)}.
\end{eqnarray}
Finally, if $p=\infty$, then $f_i=\frac12\mathbf{1}_{[-1,1]}$ and
$$
V(Z_1(\mu))\geq V(M_\infty(\mu))=2^n\int_{\R^n}^{\ast}\sup_{x=\sum_{i=1}^kc_i\theta_iu_i}  \prod_{i=1}^kf_i(\theta_i)^{c_i}\,dx
\geq 2^n\prod_{i=1}^k\left(\int_{\R}f_i\right)^{c_i}=2^n.
$$
Equality in (\ref{volZp*})  leads to equality in the reverse Brascamp-Lieb inequality, and hence $k=2n$ and
$u_1,\ldots,u_k$ form the vertices of a regular crosspolytope in $\R^n$.

The main idea in deriving a stability version of (\ref{volZ*p}) and (\ref{volZp*}) is to establish a stronger version of 
(\ref{BLstep2}) and (\ref{RBLstep2}), respectively, based on the stronger version Lemma~\ref{Ball-Barthe-stab} of (i).
In order to apply the estimate of Lemma~\ref{Ball-Barthe-stab}, we need some basic bounds on the derivatives of the transportation maps involved. These bounds are proved in  Section~\ref{secTrans}. The technical Sections~\ref{seciso} and \ref{seceveniso}
also serve as a preparation for the proof of the core statement Proposition~\ref{volZ*pstab} providing the stabiliy version
of (\ref{BLstep2}). The argument for the estimate strenghtening (\ref{RBLstep2}) is similar, and is reviewed in
Section~\ref{secLpzonoids}. This finally completes the proof of Theorem~\ref{Zpmustab}. 
The stability versions of the reverse isoperimetric inequality in the origin symmetric case 
(Theorem~\ref{vol} and Theorem\ref{BM}) and the strengthening of Theorem 1.4 for $p=\infty$ 
stated in Corollary~\ref{Zinftymustab} are proved in Section~\ref{secreverseisostab}.

The methods of this paper are very specific for our particular choice of the functions $f_i$, and
no method is known to the authors that could lead to a stability version of the Brascamp-Lieb inequality (\ref{BL}) 
or of its reverse form (\ref{RBL}) in general. However, the proof of Theorem~\ref{Zpmustab} suggests the following conjecture.

\begin{conj}
\label{BL-RBLsamefunction}
If $f$ is an even probability density function on $\R$ with variance $1$, $g(t)=\frac1{\sqrt{2\pi}}\,e^{- t^2/2}$ is the standard normal distribution, and
$\mu$ is an even isotropic measure on $S^{n-1}$ supported at $u_1,\ldots,u_k\in S^{n-1}$ with $\mu(\{u_i\})=c_i$, then
\begin{eqnarray*}
\int_{\R^n}\prod_{i=1}^kf(\langle x,u_i\rangle)^{c_i}\,dx&\leq&
\exp\left(-\gamma \min\{1,\|f-g\|_1\}^\alpha\cdot \delta_{\rm WO}(\mu,\nu_n)^\alpha\right), \\
\int_{\R^n}^{\ast}\sup_{x=\sum_{i=1}^kc_i\theta_iu_i}\prod_{i=1}^kf(\theta_i)^{c_i}\,dx&\geq&
\exp\left(\gamma \min\{1,\|f-g\|_1\}^\alpha\cdot\delta_{\rm WO}(\mu,\nu_n)^\alpha\right),
\end{eqnarray*}
where $\gamma>0$ depends on $n$ and $\alpha>0$ is an absolute constant.
\end{conj}

\section{An auxiliary analytic stability result}

To obtain a stability version of Theorem~B, we need a stability version of 
the Brascamp-Lieb inequality and its reverse form
in the special cases we use. For this we need some analytic inequalities such as 
 estimates of the derivatives of the corresponding transportation maps, which will 
be provided in Section~\ref{secTrans}. Moreover, we will use the following  
strengthened form of (i) and a basic algebraic inequality, which were both 
established in \cite[Section 4]{BoH}.

\begin{lemma}
\label{Ball-Barthe-stab}
Let $k\geq n+1$, $t_1,\ldots,t_k>0$, and let $v_1,\ldots,v_k\in \R^n$ satisfy 
$\sum_{i=1}^kv_i\otimes v_i={\rm Id}_n$. Then
$$
\det\left( \sum_{i=1}^kt_iv_i\otimes v_i\right)\geq \theta^*\ 
\prod_{i=1}^k t_i^{\langle  v_i,v_i\rangle},
$$
where
\begin{align*}
\theta^*&=  1+\frac12\sum_{1\leq i_1<\ldots<i_n\leq k}
\det[v_{i_1},\ldots,v_{i_n}]^2
\left(\frac{\sqrt{t_{i_1}\cdots  t_{i_n}}}{t_0}-1\right)^2,\\
t_0&=\sqrt{\sum_{1\leq i_1<\ldots<i_n\leq k}
t_{i_1}\cdots t_{i_n}\det[v_{i_1},\ldots,v_{i_n}]^2}.
\end{align*}
\end{lemma}

In order to estimate $\theta^*$ from below, we use the following  observation from \cite{BoH}.

\begin{lemma}
\label{xab}
If $a,b,x>0$, then
$$
(xa-1)^2+(xb-1)^2\geq \frac{(a^2-b^2)^2}{2(a^2+b^2)^2}.
$$
\end{lemma}

\section{The transportation maps}
\label{secTrans}

We note that for $p\geq 1$, we have
\begin{equation}
\label{tpint}
\int_{\R}e^{-|t|^p}\,dt=\frac2p\int_0^\infty e^{-s}s^{\frac1p-1}\,ds=2\Gamma\mbox{$(1+\frac1p)$}.
\end{equation}
Thus for $p\in[1,\infty]$, we consider the density functions
$$
\varrho_p(x)=\left\{
\begin{array}{rl}
\frac1{2\Gamma(1+\frac1p)}\,e^{-|s|^p}&\mbox{ \ if $p\in[1,\infty)$},\\[2ex]
\frac12\mathbf{1}_{[-1,1]}&\mbox{ \ if $p=\infty$.}
\end{array}
\right.
$$
In particular, $\varrho_2$ is the Gaussian density function $\pi^{-1/2}e^{-s^2}$.
In addition,  we define the transportation maps $\varphi_p,\psi_p:\R\to \R$
for $p\in[1,\infty)$, $\varphi_\infty:(-1,1)\to \R$  and $\psi_\infty:\R\to(-1,1)$
 by
\begin{eqnarray}
\label{phip}
\int_{-\infty}^t\varrho_p(s)\,ds
&=&
\int_{-\infty}^{\varphi_p(t)} \varrho_2(s)\,ds,\\
\label{psip}
\int_{-\infty}^{\psi_p(t)}\varrho_p(s)\,ds&=&\int_{-\infty}^t \varrho_2(s)\,ds.
\end{eqnarray}
Here $\varphi_p$ and $\psi_p$ are odd and inverses of each other.\\

In the following, we use that
$$
s-s^2\leq\log(1+s)\leq s\mbox{ \ if $s\geq -\frac12$},
$$
and the following properties of the $\Gamma$ function.
\begin{description}
\item{(i)} $\log \Gamma(t)$ is strictly convex for $t>0$;

\item{(ii)} $\Gamma(1)=\Gamma(2)=1$;

\item{(iii)} $\Gamma(1+\frac1{2.3})<\Gamma(1+\frac12)=\sqrt{\pi}/2$;

\item{(iv)} $\Gamma$ has a unique minimum on $(0,\infty)$ at 
$x_{\rm min}=1.4616\ldots$ with $\Gamma(x_{\rm min})=0.885603\ldots$. 
In particular, 
$\Gamma(t)>0.8856$ for $t>0$, $\Gamma$ is strictly decreasing 
on $[0,x_{\rm min}]$ and strictly increasing on $[1.5,\infty)$.
\end{description}
We deduce from (i)--(iv) that the density functions involved satisfy
\begin{equation}
\label{densityest}
\frac1{2e}\leq\varrho_p(s)<\frac1{2\cdot 0.8856}
\mbox{ \ for $p\in[1,\infty]$ and $s\in[0,1]$.}
\end{equation}
We note that $e/0.8856<3.1$, and hence 
\begin{equation}\label{leqone}
\varphi_p(s)\in [0,1)\text{ for }s\in [0,\tfrac{1}{3.1}].
\end{equation}
In fact, assuming that $\varphi_p(\tfrac{1}{3.1})\ge 1=\varphi_p(t)$, $t\in(0,\tfrac{1}{3.1}]$, we have
$$
\frac{3.1^{-1}}{2\cdot 0.8856}>\int_0^t\varrho_p(s)\, ds=\int_0^1\varrho_2(s)\, ds\ge\frac{1}{2e},
$$
a contradiction. 
Then, (\ref{densityest}) and (\ref{masstrans}) yield that 
\begin{equation}
\label{phipsiderest}
\frac1{3.1}<\varphi'_p(s),\psi'_p(s)<3.1
\mbox{ \ for $p\in[1,\infty]$ and $s\in[0,\frac1{3.1}]$.}
\end{equation}

The following simple estimate will play a crucial role in the proofs of 
Lemma~\ref{trasportation-map-phi}
and  Lemma~\ref{trasportation-map-psi}.

\begin{lemma}
\label{p2est}
For $p\in(1,3)\setminus\{2\}$ and $\nu>0$, let $f(t)=\nu t-pt^{p-1}$ 
for $t\in [0,1]$.
\begin{description}
\item{\rm (a)} If $p\in(1,2)$, $f(\tau)\leq 0$ for some $\tau\in(0,1]$ 
and $t\in(0,\tau/2]$, then 
$$
f(t)<-\frac{p(p-1)(2-p)}{2^{4-p}}\cdot t^{p-1}.
$$
\item{\rm (b)} If $p\in(2,3)$, $f(\tau)\geq 0$ for some $\tau\in(0,1]$ 
and $t\in(0,\tau/2]$, then 
$$
f(t)>\frac{p(p-1)(p-2)}{2^{4-p}}\cdot t^{p-1}.
$$
\end{description}
\end{lemma}
{\bf Remark} Naturally, the bound could be linear in $t$ with a factor depending on $\nu$, but this way the only influence of $\nu$ is on the value of $\tau$. We only use Lemma~\ref{p2est} when $1.5\leq p\leq 2.3$ and $t>c$ for a positive absolute constant $c$ anyway.\\

\bigskip

\proof Let $p\in(1,2)$.  Since $f$ is convex on $[0,\tau]$, $\tau\le 1$, $f(0)\le 0$ and $f(\tau)\le 0$, we have $f(2t)\le 0$ for $t\in [0,\tau/2]$. 
Taylor's formula yields that if $t\in(0,\tau/2]$, then there exist $\tau_1\in (0,t)$ and $\tau_2\in(t,2t)$ such that
\begin{eqnarray*}
0&\geq& \frac12\left(f(0)+f(2t)\right)=\frac12\left(f(t)-f'(t)t+\frac12\,f''(\tau_1)t^2+f(t)+f'(t)t+\frac12\,f''(\tau_2)t^2\right)\\
&=&f(t)+\frac12\,\frac{f''(\tau_1)+f''(\tau_2)}2\,t^2,
\end{eqnarray*}
where $0<\tau_i<2t\le \tau$. From $f''(\tau_i)=-p(p-1)(p-2)\tau_i^{p-3}>p(p-1)(2-p) (2t)^{p-3}$, $i=1,2$, we deduce the estimate
$$
f(t)<-\frac12\,p(p-1)(2-p) (2t)^{p-3}\cdot t^2=
-\frac{p(p-1)(2-p)}{2^{4-p}}\cdot t^{p-1}.
$$

If $p\in(2,3)$, then  $f(t)=\nu t-pt^{p-1}$ is concave on $[0,\tau]$, and a similar argument yields (b). 
\proofbox

\begin{lemma}
\label{trasportation-map-phi}
Let $p\in[1,\infty]\setminus\{2\}$ and 
 $t\in(0,\frac18)$. Then
\begin{eqnarray}
\varphi_p''(t)&<&-\frac{2-p}{48}\cdot t\mbox{ \ \ \ if\,  $p\in[1,2)$},\label{eqa1}\\
\varphi_p''(t)&>&\frac{p-2}{5}\cdot t^{1.3}\mbox{ \ \ \ if\,  $p\in(2,3]$},\label{eqa2}\\
\varphi_p''(t)&>&0.2\cdot t^{1.3}\mbox{ \ \ \ if\, $p\in(3,\infty]$}.\label{eqa3}
\end{eqnarray}
\end{lemma}
\proof For brevity of notation, let  $\varphi=\varphi_p$.
 We have $\varphi(0)=0$ as $\varphi$ is odd. Since $\varphi$ is strictly increasing,  $\varphi(t)>0$ if $t>0$.

Let $p\in[1,\infty)\setminus \{2\}$. For $t>0$, differentiating (\ref{phip}) yields the formula
$$
\frac{e^{-t^p}}{2\Gamma(1+\frac1p)}=\frac{e^{-\varphi(t)^2}\varphi'(t)}{2\Gamma(1+\frac12)},
$$
and by differentiating again, we obtain
$$
\frac{-p\Gamma(1+\frac12)}{\Gamma(1+\frac1p)}\cdot e^{-t^p}t^{p-1}=-2 e^{-\varphi(t)^2}\varphi(t)\varphi'(t)^2
+e^{-\varphi(t)^2}\varphi''(t).
$$
In particular,
\begin{eqnarray}
\label{phider}
\varphi'(t)&=&\frac{\Gamma(1+\frac12)}{\Gamma(1+\frac1p)}\,e^{\varphi(t)^2-t^p},\\
\label{phi2der}
\varphi''(t)&=&(2\varphi(t)\varphi'(t)-pt^{p-1})\varphi'(t).
\end{eqnarray}

In the following argument, we use the value
$$
t_p=(2/p)^{\frac1{p-2}} \mbox{ \ \ \ for $p\in[1,\infty)\setminus\{2\}$}.
$$
The function $p\mapsto t_p$ is continuously extended to $p=2$ by $t_2=e^{-1/2}$, and then this function 
is increasing on $[1,\infty)$. In particular, $t_p\ge 1/2$ for $p\in[1,\infty)$.

Moreover, we apply the fact that
\begin{equation}
\label{plnt}
\mbox{for given $t\in(0,1/e)$,  $p\mapsto p t^{p-1}$ is a decreasing function of $p\geq 1$}.
\end{equation}
First, we show that for $1\leq p<2$ and $t\in(0,1/4)$, we have $\varphi''(t)<-\frac{2-p}{48}\cdot t$, which proves \eqref{eqa1}. 

In this case, $\varphi'(0)< 1$ by (\ref{phider}), (i), (ii) and (iv). Since $\varphi'$ is continuous, 
there exists a largest $s_p\in(0,\infty]$ such that $\varphi'(t)< 1$ if $0<t<s_p$. Thus, if $t\in (0,s_p)$, then $\varphi(t)< t$, and in turn
(\ref{phi2der}) yields that
$$
\varphi''(t)=(2\varphi(t)\varphi'(t)-pt^{p-1})\varphi'(t)<(2t-pt^{p-1})\varphi'(t).
$$
For $1\le p<2$ and $t\in[0,t_p]$,  we have $2t-pt^{p-1}\leq 0$.  
In particular, $\varphi'(t)$ is monotone decreasing on $(0,\min\{s_p,t_p\})$, which in turn implies that $s_p\geq t_p$. 
We deduce from (\ref{phipsiderest})  that
\begin{equation}
\label{phi2dercase1}
\varphi''(t)<\frac{2t-pt^{p-1}}{3.1} \mbox{ \ \ \ \ \ for $t\in(0,\frac{1}{3.1})$.}
\end{equation}

Now we distinguish two cases. If $1.5\leq p<2$, then we deduce from
(\ref{phi2dercase1}) and  Lemma~\ref{p2est} (a) that
\begin{equation}
\label{phi2dercase10}
\varphi''(t)<-\frac{p(p-1)(2-p)}{3.1\cdot 2^{4-p}}\cdot t^{p-1}<
-\frac{\frac34(2-p)}{3.1\cdot 2^{2.5}}\cdot t<-\frac{2-p}{24}\cdot t
\mbox{ \ \ \ \ \ for $t\in(0,\frac1{4})$.}
\end{equation}
If $1\leq p\leq 1.5$, then when estimating the right-hand side of (\ref{phi2dercase1}) for a given $t\in(0,\frac1{4})$, 
we may assume that $p=1.5$ according to (\ref{plnt}). In other words, using Lemma~\ref{p2est} (a), inequality (\ref{phi2dercase10}) yields that if $1\leq p\leq 1.5$ and 
$t\in(0,\frac1{4})$, then
$$
\varphi''(t)<\frac{2t-pt^{p-1}}{3.1}\leq \frac{2t-1.5 t^{0.5}}{3.1}\leq 
-\frac{2-1.5}{24}\cdot t\leq-\frac{2-p}{48}\cdot t.
$$

Second, if $2<p\leq 2.3$ and $t\in(0,\frac14)$, then we show that $\varphi''(t)>\frac{p-2}{2}\cdot t^{1.3}$.

In this case, $\varphi'(0)> 1$ by (\ref{phider}), (i), (iii)  and (iv). Since $\varphi'$ is continuous, there exists a largest $s_p\in(0,\infty]$ such that $\varphi'(t)> 1$ if $0<t<s_p$. Thus if $t\in (0,s_p)$, then $\varphi(t)> t$, and in turn
(\ref{phi2der}) yields that
$$
\varphi''(t)=(2\varphi(t)\varphi'(t)-pt^{p-1})\varphi'(t)>(2t-pt^{p-1})\varphi'(t).
$$
For $p>2$ and $t\in[0,t_p]$, we have  $2t-pt^{p-1}\geq 0$. 
In particular, $\varphi'(t)$ is monotone increasing on $(0,\min\{s_p,t_p\})$, which, in turn, implies that $s_p\geq t_p$. We deduce  that
\begin{equation}
\label{phi2dercase2}
\varphi''(t)>2t-pt^{p-1} \mbox{ \ \ \ \ \ if $t\in(0,\frac12)$.}
\end{equation}
We deduce from (\ref{phi2dercase2}) and  Lemma~\ref{p2est} (b) that
$$
\varphi''(t)>\frac{p(p-1)(p-2)}{2^{4-p}}\cdot t^{p-1}>
\frac{2(p-2)}{2^{2}}\cdot t^{1.3}=\frac{p-2}{2}\cdot t^{1.3}
\mbox{ \ \ \ \ \ if $t\in(0,\frac14)$.}
$$

 If $p\geq 2.3$ and $t\in(0,\frac18)$, then $\varphi''(t)>0.2\cdot t^{1.3}$, which completes the proof of \eqref{eqa2}. 

In this case, $\varphi'(0)> \sqrt{\pi}/2$ by (\ref{phider}), (i)--(iv). Since $\varphi'$ is continuous, there exists largest $s_p\in(0,\frac14]$ such that $\varphi'(t)> \sqrt{\pi}/2$ if $0<t<s_p$. Thus if $t\in (0,s_p]$, then 
$\varphi(t)> (\sqrt{\pi}/2)\cdot t$. From \eqref{plnt} we see that 
$$
2\varphi(t)\varphi'(t)-pt^{p-1}\ge \frac{\pi}2\,t-pt^{p-1}\ge \frac{\pi}2\,t-2.3t^{1.3}\ge 0 
$$
for $0<t\le s_p\le1/4$. Hence 
(\ref{phi2der}) yields that
$$
\varphi''(t)=(2\varphi(t)\varphi'(t)-pt^{p-1})\varphi'(t)>\left(\frac{\pi}2\,t-2.3t^{1.3}\right)\cdot \frac{\sqrt{\pi}}2
$$
for $t\in(0,s_p]$. In particular, we conclude that $s_p=\frac14$, and hence
Lemma~\ref{p2est} (b) yields that
$$
\varphi''(t)>\frac{(\sqrt{\pi}/2)\cdot 2.3\cdot 1.3\cdot 0.3}{2^{1.7}}\cdot t^{1.3}>
0.2\cdot t^{1.3}
\mbox{ \ \ \ \ \ for $t\in(0,\frac18)$.}
$$

If $p=\infty$ and $t>0$, then $\varphi''(t)>t$, which completes the proof of \eqref{eqa3}.
Differentiating (\ref{phip}) we deduce for $t\in (-1,1)$ that
\begin{eqnarray}
\label{phiinftyder}
\varphi'(t)&=&\Gamma\left(1+\frac12\right)e^{\varphi(t)^2}=\frac{\sqrt{\pi}}2\,e^{\varphi(t)^2}, \\
\label{phiinftyder2}
\varphi''(t)&=&2\varphi(t)\varphi'(t)^2.
\end{eqnarray}
As $\varphi(t)>0$ for $t>0$, we have $\varphi''(t)\geq 0$ by (\ref{phiinftyder2}), and hence $\varphi'(t)$ is monotone increasing for $t\geq 0$. Therefore $\varphi'(t)\geq \varphi'(0)=\sqrt{\pi}/2$ by (\ref{phiinftyder}), which, in turn,  again by 
(\ref{phiinftyder2}) yields that
$$
\varphi''(t)\geq 2\left(\frac{\sqrt{\pi}}2\right)^3t>t\mbox{ \ \ \ \ for $t\in (0,1)$.}
$$
Thus we have proved all estimates of Lemma~\ref{trasportation-map-phi} for $\varphi^{\prime\prime}$. 
\proofbox

\begin{lemma}
\label{trasportation-map-psi}
Let $p\in[1,\infty]\setminus\{2\}$.
For $t\in(0,\frac1{10})$, we have
\begin{eqnarray}
\psi_p''(t)&>&\frac{2-p}{16}\cdot t\mbox{ \ \ \ if $p\in[1,2)$},\label{eqb1}\\
\psi_p''(t)&<&-\frac{p-2}{11}\cdot t^{1.3}\mbox{ \ \ \ if $p\in(2,3]$},\label{eqb2}\\
\psi_p''(t)&<&- \frac{1}{11}\cdot t^{1.3}\mbox{ \ \ \ if $p\in(3,\infty]$.}\label{eqb3}
\end{eqnarray}
\end{lemma}
\proof To simplify notation, let  $\psi=\psi_p$.
 We have  $\psi(0)=0$ as $\psi$ is odd. Therefore $\psi(t)>0$ if $t>0$.
Turning to $\psi''$, we only sketch the main steps. In this case, 
differentiating (\ref{psip}) yields the formulas
\begin{eqnarray}
\label{psider}
\psi'(t)&=&\frac{\Gamma(1+\frac1p)}{\Gamma(1+\frac12)}\,e^{\psi(t)^p-t^2},\nonumber\\
\label{psi2der}
\psi''(t)&=&(p\psi(t)^{p-1}\psi'(t)-2t)\psi'(t).
\end{eqnarray}

First, for $1\leq p<2$ and $t\in(0,\frac1{8})$ we show that $\psi''(t)>\frac{2-p}{16}\cdot t$, which proves \eqref{eqb1}. 

If $p\in[1,2)$, then $\psi'(0)>1$ by (i), (ii) and (iv). Arguments similar to those in the proof of 
 Lemma~\ref{trasportation-map-phi} yield
\begin{equation}
\label{psi2dercase1}
\psi''(t)=(p\psi(t)^{p-1}\psi'(t)-2t)\psi'(t)>pt^{p-1}-2t \mbox{ \ \ \ \ \ for $t\in(0,\frac12)$.}
\end{equation}
 If $1.5\leq p<2$, then we deduce from
(\ref{psi2dercase1}) and  Lemma~\ref{p2est} (a) that
$$
\psi''(t)>\frac{p(p-1)(2-p)}{2^{4-p}}\cdot t^{p-1}>
\frac{\frac34(2-p)}{2^{2.5}}\cdot t>\frac{2-p}{8}\cdot t
\mbox{ \ \ \ \ \ for $t\in(0,\frac1{8})$.}
$$

If $1\leq p\leq 1.5$, then when estimating the right-hand side of (\ref{psi2dercase1}) for a given $t\in(0,\frac1{e})$, 
we may assume that $p=1.5$ according to (\ref{plnt}). In other words, (\ref{psi2dercase1}) yields that if $1\leq p\leq 1.5$ and 
$t\in(0,\frac1{e})$, then
\begin{equation}
\label{psi2dercase11}
\psi''(t)>pt^{p-1}-2t\geq 1.5 t^{0.5}-2t\geq 
\frac{2-1.5}{8}\cdot t\geq \frac{2-p}{16}\cdot t.
\end{equation}

Next, for $2<p\leq 2.3$ and $t\in(0,\frac14)$, we prove that $\psi''(t)<-\frac{p-2}{3}\cdot t^{1.3}$.

If $p\in(2,2.3]$, then $\psi'(0)<1$ by (i)--(iv), and arguments similar to the ones used in the proof of 
  Lemma~\ref{trasportation-map-phi} yield
$$
\psi''(t)=(p\psi(t)^{p-1}\psi'(t)-2t)\psi'(t)<-(2t-pt^{p-1})\psi'(t)<-\frac{2t-pt^{p-1}}{3.1}<0
\mbox{ \ \ \ \ \ for $t\in(0,\frac1{3.1})$.}
$$
We deduce from  Lemma~\ref{p2est} (b) that
$$
\psi''(t)<-\frac{p(p-1)(p-2)}{3.1\cdot 2^{4-p}}\cdot t^{p-1}<
-\frac{2(p-2)}{3.1\cdot 2^{2}}\cdot t^{1.3}<-\frac{p-2}{7}\cdot t^{1.3}
\mbox{ \ \ \ \ \ for $t\in(0,\frac18)$.}
$$

Let $p\geq 2.3$ and $t\in(0,\frac1{10})$. We now show that $\psi''(t)<- t^{1.3}/11$, which completes the proof of \eqref{eqb2}.

In this case, $\psi'(0)< 2/\sqrt{\pi}$ by  (i)--(iv).
There exists a maximal $s_p\in (0,\frac15]$ such that if $t\in(0,s_p)$, then
 $\psi'(t)<2/\sqrt{\pi}$. Thus if $t\in (0,s_p]$, then 
$\psi(t)< (2/\sqrt{\pi})\cdot t$, and, in turn,
(\ref{psi2der}) yields that
\begin{equation}
\label{psi2dercase3zero}
\psi''(t)=(p\psi(t)^{p-1}\psi'(t)-2t)\psi'(t)<\left(\left(\frac2{\sqrt{\pi}}\right)^ppt^{p-1}-2t\right)\psi'(t).
\end{equation}
Given $t\in(0,\frac12]$,
$$
\frac{d}{dp}\,\log \left[\left(\frac2{\sqrt{\pi}}\right)^ppt^{p-1}\right]=
\frac1p+\log \frac{2t}{\sqrt{\pi}} <0 \mbox{ \ \ for $p\in(2,\infty)$},
$$
and hence (\ref{psi2dercase3zero}) yields that  if $t\in (0,s_p]$, then 
\begin{eqnarray}
\label{psi2dercase3one}
\psi''(t)&=&(p\psi(t)^{p-1}\psi'(t)-2t)\psi'(t)\nonumber\\
&<&\left(\left(\frac2{\sqrt{\pi}}\right)^{2.3}2.3t^{1.3}-2t\right)\psi'(t)=
f(t)\left(\frac2{\sqrt{\pi}}\right)^{2.3}\psi'(t)
\end{eqnarray}
where
$$
f(t)=2.3t^{1.3}-2\left(\frac{\sqrt{\pi}}2\right)^{-2.3}t.
$$
Here $f(\frac15)<0$, thus with $\tau=\frac15$,  Lemma~\ref{p2est} (b) yields that
$$
f(t)<- \frac{2.3\cdot 1.3\cdot 0.3}{2^{1.7}}\cdot t^{1.3}<-0.27\cdot t^{1.3}
\mbox{ \ \ \ \ \ for $t\in(0,\frac1{10})$.}
$$
We conclude from (\ref{phipsiderest}) and (\ref{psi2dercase3one}) that
$$
\psi''(t)<-\frac{(\frac2{\sqrt{\pi}})^{2.3}\cdot 0.27\cdot t^{1.3}}{3.1}<-\frac{t^{1.3}}{11}
\mbox{ \ \ \ \ \ for $t\in(0,\frac1{10})$.}
$$

Finally, for $p=\infty$ and $t\in(0,\frac1{3.1})$, we show $\psi''(t)<-\frac{2}{3.1}\cdot t$, which 
completes the proof of \eqref{eqb3}.

Differentiating (\ref{psip}) we deduce that if $t>0$, then
\begin{eqnarray*}
\psi'(t)&=&\frac1{\Gamma\left(1+\frac12\right)}e^{-t^2}=\frac2{\sqrt{\pi}}\,e^{-t^2}, \\
\psi''(t)&=&-2t\psi'(t).
\end{eqnarray*}
We conclude from (\ref{phipsiderest})  that $\psi''(t)<-\frac{2t}{3.1}$ for $t\in(0,\frac1{3.1})$.

In summary, we have established all estimates of Lemma~\ref{trasportation-map-psi} for $\psi''$. 
\proofbox

\section{Basic estimates on isotropic measures}
\label{seciso}

The main result of this section is Lemma~\ref{nsmallballs}. It states that for any isotropic measure $\mu$ on $S^{n-1}$, there exist 
spherical caps $X_1,\ldots,X_n\subseteq S^{n-1}$ whose $\mu$-measure is bounded from below and which have the additional property 
that for any vectors $w_i\in X_i$, $i\in\{1,\ldots,n\}$, also the determinant 
$|\det[w_1,\ldots,w_n]|$ is bounded from below.

For $\alpha\in(0,\frac{\pi}2]$ and $v\in S^{n-1}$, we consider the closed and open spherical caps
\begin{eqnarray*}
\Omega(v,\alpha)&=&\{u\in S^{n-1}:\,\langle u,v\rangle\geq \cos\alpha\},\\
\widetilde{\Omega}(v,\alpha)&=&\{u\in S^{n-1}:\,\langle u,v\rangle> \cos\alpha\}.
\end{eqnarray*}

\begin{claim}
\label{isotropiccap}
If $\mu$ is an isotropic measure on $S^{n-1}$, $v\in S^{n-1}$, and 
$\alpha\in(0,\frac{\pi}2)$, then
$$
\mu\left(\widetilde{\Omega}(v,\alpha)\right) +
\mu\left(\widetilde{\Omega}(-v,\alpha)\right)\geq  1-n\cos^2\alpha.
$$
\end{claim}
\proof For given $v\in S^{n-1}$ and 
$\alpha\in(0,\frac{\pi}2)$, let $X=\{u\in S^{n-1}:\,|\langle u,v\rangle|\leq \cos\alpha\}$. Since $\mu$ is isotropic, 
we have $\mu(X)\leq n$, and
\begin{eqnarray*}
1&=&\langle v,v\rangle=\int_{S^{n-1}}\langle u,v\rangle^2\,d\mu(u)=
\int_{\widetilde{\Omega}(v,\alpha)\cup \widetilde{\Omega}(-v,\alpha)}\langle u,v\rangle^2\,d\mu(u)+
\int_{X}\langle u,v\rangle^2\,d\mu(u)\\
&\leq &\mu\left(\widetilde{\Omega}(v,\alpha)\cup \widetilde{\Omega}(-v,\alpha)\right)+n\cos^2\alpha.
\mbox{ \ }\proofbox
\end{eqnarray*}

Observe that if  $\cos\alpha\ge 1/\sqrt{n}$ in the preceding claim, then the conclusion holds trivially. 

The next claim follows from a standard argument but we are not aware of any reference. 

\begin{claim}
\label{bigsmallcap}
If $\mu$ is a Borel measure on $S^{n-1}$, $p\in S^{n-1}$, and 
$0<\beta<\alpha<\frac{\pi}2$, then there exists a point $v\in \Omega(p,\alpha)$ such that
$$
\mu\left(\Omega(p,\alpha)\cap \Omega(v,\beta)\right)\ge \mu(\Omega(p,\alpha))
\cdot \frac{\sin^{n-1}\beta}{\sqrt{2\pi n}};
$$
if $\mu(\Omega(p,\alpha))>0$, then $v\in \Omega(p,\alpha)$ can be chosen such that the inequality is strict.
\end{claim}
\proof We define the Borel measure $\bar{\mu}$ on $S^{n-1}$ by $\bar{\mu}(X)=\mu(X\cap \Omega(p,\alpha))$ 
for  Borel sets $X\subseteq S^{n-1}$. 
Let $\nu$ be the Haar probability measure on ${\rm SO}(n)$. Hence, if $X\subseteq S^{n-1}$ is a Borel set and $u\in S^{n-1}$, then
$$
\nu(\{g\in {\rm SO}(n):\,gu\in X\}=\frac{{\cal H}^{n-1}(X)}{{\cal H}^{n-1}(S^{n-1})},
$$
where $\mathcal{H}^{n-1}$ denotes the $(n-1)$-dimensional Hausdorff measure (its restriction to Borel subsets of $S^{n-1}$ 
equals spherical Lebesgue measure). 
We deduce that
\begin{eqnarray*}
\mu\left(\Omega(p,\alpha)\right)\cdot 
\frac{{\cal H}^{n-1}(\Omega(p,\beta))}{{\cal H}^{n-1}(S^{n-1})}
&=&\bar{\mu}\left(S^{n-1}\right)\cdot 
\frac{{\cal H}_{n-1}(\Omega(p,\beta))}{{\cal H}_{n-1}(S^{n-1})}\\
&=&
\int_{S^{n-1}}\int_{{\rm SO}(n)}\mathbf{1}_{\Omega(p,\beta)}(gu)\,d\nu(g)\,d\bar{\mu}(u)\\
&=&
\int_{{\rm SO}(n)}\int_{S^{n-1}}\mathbf{1}_{\Omega(p,\beta)}(gu)\,d\bar{\mu}(u)\,d\nu(g)\\
&=&\int_{{\rm SO}(n)}\bar{\mu}(\Omega(g^{-1}p,\beta))\,d\nu(g)\\
&=&
\int_{{\rm SO}(n)}\mu(\Omega(p,\alpha)\cap\Omega(g^{-1}p,\beta))\,d\nu(g).
\end{eqnarray*}
Hence there exists some $v_0\in S^{n-1}$ such that
$$
\mu\left(\Omega(p,\alpha)\cap \Omega(v_0,\beta)\right)
\geq \mu\left(\Omega(p,\alpha)\right)\cdot 
\frac{{\cal H}^{n-1}(\Omega(p,\beta))}{{\cal H}^{n-1}(S^{n-1})}.
$$
To finish the proof, we can assume that $\mu(\Omega(p,\alpha))>0$. 
Finally, if $v\in \Omega(p,\alpha)$ is the closest point to $v_0$, then
$$
\Omega(p,\alpha)\cap \Omega(v_0,\beta)\subseteq
\Omega(p,\alpha)\cap \Omega(v,\beta). 
$$
To conclude the proof, we use that 
${\cal H}^{n-1}(\Omega(p,\beta))>\kappa_{n-1}\sin^{n-1}\beta$,
${\cal H}^{n-1}(S^{n-1})=n\kappa_n$, where $\kappa_i$ denotes 
the volume of the $i$-dimensional unit ball, 
and the basic inequality 
$\frac{\kappa_{n-1}}{n\kappa_n}>\frac1{\sqrt{2\pi n}}$, which follows 
from (i); see \cite[p.~564, l.~2]{Wendel48}. 
\proofbox

\begin{claim}
\label{largedeterminant}
If $b_1,\ldots,b_n\in S^{n-1}$, and $s_1,\ldots,s_n\in\R^n$ satisfy 
$\|s_i\|\leq|\det[b_1,\ldots,b_n]|/4n$, then
$$
|\det[b_1+s_1,\ldots,b_n+s_n]|\geq |\det[b_1,\ldots,b_n]|/2.
$$
\end{claim}
\proof Let $D=|\det[b_1,\ldots,b_n]|/4n$.
Since for 
any $r_1,\ldots,r_n\in\R^n$ we have
$$
|\det[r_1,\ldots,r_n]|\leq \|r_1\|\cdots  \|r_n\|,
$$
we deduce from the linearity of the determinant
and $e^t\le 1+2t$ for $t\in [0,1]$ that
\begin{eqnarray*}
|\det[b_1+s_1,\ldots,b_n+s_n]|&\geq &
|\det[b_1,\ldots,b_n]|-\sum_{i=1}^n{n\choose i}D^i\\
&=&4nD-(1+D)^n+1\\
&\geq &4nD-e^{nD}+1\\
&\geq & 4nD-2nD\geq 2nD=|\det[b_1,\ldots,b_n]|/2. \mbox{ \ \ }\proofbox
\end{eqnarray*}

Lemma~\ref{nsmallballs} can be considered as a measure theoretic version of the Dvoretzky-Rogers lemma
(see A. Dvoretzky, C. A. Rogers \cite{DvR50}, S. Brazitikos, A. Giannopoulos, P. Valettas, 
B.-H. Vritsiou \cite{BGVV16}, and for a non-symmetric version, M. Naszodi \cite{Nas16}).

\begin{lemma}
\label{nsmallballs} 
Let $\beta=2^{-(n+1)}n^{-(n+1)/2}$.
If $\mu$ is an isotropic measure on $S^{n-1}$, then there exist $v_1,\ldots,v_n\in S^{n-1}$
such that $\mu(\Omega(v_i,\beta))\geq \beta^n$, for $i=1,\ldots,n$, and such that if
$w_i\in \Omega(v_i,\beta)$, for $i\in\{1,\ldots,n\}$, then
$|\det[w_1,\ldots,w_n]|\geq2n\beta$.
\end{lemma}
\proof Let $\alpha_n\in(0,\frac{\pi}2)$ satisfy $\cos\alpha_n=\frac1{2\sqrt{n}}$.
First, we will construct $v_i,p_i\in S^{n-1}$ by induction on $i\in \{1,\ldots,n\}$ in such a way that
\begin{eqnarray}
\label{vipi1}
\mu(\Omega(v_i,\beta))&\geq&\beta^n,\\
\label{vipi2}
\mu(\Omega(p_i,\alpha_n))&\geq &3/8,\\
\label{vipi3}
v_i&\in&\Omega(p_i,\alpha_n),\\
\label{vipi4}
\langle p_i,v_j\rangle&=&0\mbox{ \ for $1\leq j<i\leq n$}.
\end{eqnarray}

For this, let $p\in S^{n-1}$. According to Claim~\ref{isotropiccap},  we can choose $p_1\in\{p,-p\}$ such that
$$
\mu(\Omega(p_1,\alpha_n))\geq \frac{1-n\cos^2\alpha_n}2=\frac38.
$$
Thus, since $\beta<1<\alpha_n$, Claim~\ref{bigsmallcap} yields the existence of
a point $v_1\in\Omega(p_1,\alpha_n)$ satisfying (\ref{vipi1}). 

If $i\geq 2$, and $v_j,p_j$ are known for $j=1,\ldots,i-1$, then we  choose $p'_i\in S^{n-1}$ satisfying 
 (\ref{vipi4}). Again, Claim~\ref{isotropiccap} provides $p_i\in\{p'_i,-p'_i\}$ satisfying (\ref{vipi2}). In addition,
 a point $v_i\in\Omega(p_i,\alpha_n)$ satisfying (\ref{vipi1}) is provided by 
Claim~\ref{bigsmallcap}.

We deduce from  (\ref{vipi3}) that if $i\in \{1,\ldots, n\}$, then 
$\langle p_i,v_i\rangle\geq\frac1{2\sqrt{n}}$.  Combined with  (\ref{vipi4}), for $i\in \{2,\ldots, n\}$ this
yields that
$$
{\rm dist}\,(v_i,{\rm aff}\,\{v_1,\ldots,v_{i-1}\})\geq \frac1{2\sqrt{n}}.
$$
In particular,
$$
|\det[v_1,\ldots,v_n]|\geq 2^{-(n-1)}n^{-(n-1)/2}=4n\beta.
$$
Next let $w_i\in \Omega(v_i,\beta)$ for $i=1,\ldots,n$, and hence
$\|s_i\|<\beta$ for $s_i=w_i-v_i$ and $i=1,\ldots, n$. Therefore Claim~\ref{largedeterminant}
implies the lemma.
\proofbox

The following Lemma~\ref{ntinyballs} uses the notation of Lemma~\ref{nsmallballs}. 

\begin{lemma}
\label{ntinyballs} 
For an  isotropic measure $\mu$ on $S^{n-1}$, let $v_1,\ldots,v_n\in S^{n-1}$ and $\beta$ be as in
Lemma~\ref{nsmallballs}. For every $i\in\{1,\ldots,n\}$ and $\eta\in(0,\beta)$, 
\begin{description}
\item{\rm (i)}  there exists $q_i\in\Omega(v_i,\beta)$ such that
$$
\mu(\Omega(v_i,\beta)\cap \Omega(q_i,\eta))\geq \frac{\beta^n}{4n},
$$
\item{\rm (ii)} or there exist $\Psi_1,\Psi_2\subseteq \Omega(v_i,\beta)$ such that
\begin{eqnarray*}
\mu(\Psi_j)&\geq&\frac{\beta^n}{4n} \mbox{ \ for $j=1,2$},\\
\|a_1-a_2\|&\geq &\frac{\eta}{\sqrt{n}}
\mbox{ \ for $a_1\in\Psi_1$ and $a_2\in \Psi_2$}.
\end{eqnarray*}
\end{description}
The  points $q_1,q_2$ and the sets $\Psi_1,\Psi_2$ can be chosen independently of $\eta\in(0,\beta)$. 
\end{lemma}
\proof Let $i\in\{1,\ldots,n\}$ and $\eta\in(0,\beta)$ be fixed. 

If there exists $q_i\in \Omega(v_i,\beta)$ such that 
$\mu(\{q_i\})\geq \frac{\beta^n}{4n}$, then (i) is satisfied. Therefore we assume that
\begin{equation}
\label{noheavypoint}
\mu(\{q\})< \frac{\beta^n}{4n} \mbox{ \ for all $q\in \Omega(v_i,\beta)$}.
\end{equation}

We choose an orthonormal basis $w_1,\ldots,w_{n-1}$ for $v_i^\bot$. 
It follows from (\ref{noheavypoint}) that there exist  
 $-1<s_j\leq t_j<1$ for $j=1,\ldots,n-1$ such that
\begin{align*}
\mu\left(\{x\in \Omega(v_i,\beta):\,\langle w_j,x\rangle<s_j \} \right)
&\leq\frac{\beta^n}{4n}\leq 
\mu\left(\{x\in \Omega(v_i,\beta):\,\langle w_j,x\rangle\leq s_j \} \right)\\
\mu\left(\{x\in \Omega(v_i,\beta):\,\langle w_j,x\rangle>t_j \} \right)
&\leq \frac{\beta^n}{4n}\leq 
\mu\left(\{x\in \Omega(v_i,\beta):\,\langle w_j,x\rangle\geq t_j \} \right).
\end{align*}
We may assume that $t_1-s_1\geq \ldots\geq t_{n-1}-s_{n-1}$, and 
we define
$\Psi_1=\{x\in \Omega(v_i,\beta):\,\langle w_1,x\rangle\leq s_1 \}$ and 
$\Psi_2=\{x\in \Omega(v_i,\beta):\,\langle w_1,x\rangle\geq t_1\}$. In addition, 
let $q_i\in\Omega(v_i,\beta)$ be such that
$\langle q_i,w_j\rangle=(s_j+t_j)/2$ for $j=1,\ldots,n-1$, and let 
$$
\Psi=\{x\in \Omega(v_i,\beta):s_j\leq  \langle w_j,x\rangle\leq t_j,\;j=1,\ldots,n-1\}.
$$

If  $t_1-s_1\geq \eta/\sqrt{n}$, then $\Psi_1$ and $\Psi_2$ satisfy (ii). 
Finally, we assume that $t_1-s_1< \eta/\sqrt{n}$, and hence
$t_j-s_j< \eta/\sqrt{n}$ for $j=1,\ldots,n-1$. 
On the one hand,
$$
\mu(\Psi)\geq\mu(\Omega(v_i,\beta))-2n\cdot \frac{\beta^n}{4n}\geq\frac{\beta^n}2.
$$
On the other hand, $\|z-(q_i|v_i^\bot)\|\leq \eta/2$ for $z\in \Psi|v_i^\bot$. Since
$\langle u,v_i\rangle>1/2$ for $u\in \Omega(v_i,\beta)$, we deduce that
$\Psi\subseteq \Omega(q_i,\eta)$. In turn, we conlude (i).
\proofbox

\section{Even isotropic measures and the cross measure}
\label{seceveniso}

As a consequence of Claim~\ref{isotropiccap}, we  estimate  the Wasserstein distance.

\begin{lemma}
\label{Hausdorffcrosserror}
Let  $\mu$ be an even isotropic measure, and let $\nu$ be a cross measure on $S^{n-1}$
with ${\rm supp}\,\nu=\{\pm w_1,\ldots,\pm w_n\}$. If 
$\delta\in[0, \frac{\pi}{4})$ and $\omega\in[0,1)$ are such that 
$$
\mu\left(S^{n-1}\backslash\bigcup_{i=1}^n(\Omega(w_i,\delta)\cup\Omega(-w_i,\delta))\right)
\leq \omega,
$$
 then
$$
\delta_W(\mu,\nu)\leq 2n \delta+2\pi n^2\omega.
$$
\end{lemma}
\proof
We write $w_{i+n}=-w_i$ for $i=1,\ldots,n$. Since $\widetilde{\Omega}\left(w_i,\frac{\pi}2-\delta\right)$
is disjoint from $\Omega(w_j,\delta)$ for $i\neq j$, it follows from 
Claim~\ref{isotropiccap} that for each $i=1,\ldots,n$, we have
\begin{eqnarray*}
\mu\left(\Omega(w_i,\delta)\cup \Omega(-w_i,\delta)\right)&\geq&
\mu\left(\widetilde{\Omega}\left(w_i,\frac{\pi}2-\delta\right)\cup 
\widetilde{\Omega}\left(-w_i,\frac{\pi}2-\delta\right)\right)-\omega\\
&>&1-n\sin^2\delta-\omega>1-n\delta^2-\omega.
\end{eqnarray*}
Since $\mu$ is even, we get
$$
\mu\left(\Omega(w_i,\delta)\right)-\frac{1}{2}\ge - \frac{n\delta^2+\omega}{2}.
$$
Since $\mu(S^{n-1})=n$, $\mu$ is even, and $\delta<\pi/4$ we deduce for $i=1,\ldots,n$ that
\begin{eqnarray*}
n&\ge&2\mu\left(\Omega(w_i,\delta)\right)+\sum_{j:j\notin\{ i, i+n\}}\mu\left(\Omega(w_j,\delta)\cup \Omega(-w_j,\delta)\right)+0\\
&\ge& 2\mu\left(\Omega(w_i,\delta)\right)+(n-1)(1-n\delta^2-\omega),
\end{eqnarray*}
and hence
$$
\mu(\Omega(w_i,\delta))\leq\frac12\left(n-(n-1)(1-n\delta^2-\omega)\right)\leq
\frac{1+n^2\delta^2+n\omega}2,
$$
for $i=1,\ldots,2n$. Thus, for $i=1,\ldots,2n$ we get
$$
\left|\mu(\Omega(w_i,\delta))-\frac{1}{2}\right|\le  \frac{n^2\delta^2+n\omega}{2}.
$$
For $f\in {\rm Lip}_1(S^{n-1})$, we may assume that $f(w_1)=0$, since $\mu(S^{n-1})=\nu(S^{n-1})=n$, and hence $|f(u)|\leq \pi$
for $u\in S^{n-1}$. Therefore 
\begin{eqnarray*}
&&\int_{S^{n-1}}f\,d\mu-\int_{S^{n-1}}f\,d\nu\\
&&\qquad =
\sum_{i=1}^{2n}
\left(\int_{\Omega(w_i,\delta)}(f(u)-f(w_i))\,d\mu(u)+\int_{\Omega(w_i,\delta)}f(w_i)\,d\mu(u)-\frac{f(w_i)}2\right)\\
&&\qquad\qquad +\int_{S^{n-1}\backslash(\cup_{i=1}^{2n}\Omega(w_i,\delta))}f(u)\,d\mu(u)\\
& &\qquad \le 2n\left(\delta\cdot\frac{1+n^2\delta^2+n\omega}2+\pi\cdot \frac{n^2\delta^2+n\omega}2\right)+\pi\omega\\
&&\qquad \leq 2n\delta+2\pi n^2\omega, 
\end{eqnarray*}
which yields the assertion.
\proofbox

We deduce the following estimate for the Wasserstein distance.

\begin{coro}
\label{Hausdorffcross}
If  $\mu$ is an even isotropic measure, and $\nu$ is a cross measure on $S^{n-1}$, and 
$\delta_H({\rm supp}\,\mu,{\rm supp}\,\nu)<\pi/4$, then
$$
\delta_W(\mu,\nu)\leq 2n \delta_H({\rm supp}\,\mu,{\rm supp}\,\nu).
$$
\end{coro}

Finally, we consider the stability of optimal symmetric coverings of $S^{n-1}$ by $2n$ 
congruent spherical caps, where a symmetric covering is an arrangement invariant under 
the antipodal map. It is a well-known conjecture that  in an optimal covering of $S^{n-1}$ 
by $2n$ congruent spherical caps, the spherical centers of the caps are vertices of a regular crosspolytope 
(see, say, L. Fejes T\'oth \cite{LFT72}). This conjecture has been verified by L. Fejes T\'oth 
\cite{LFT72} for $n\leq 3$, and by 
L. Dalla, D. G Larman, P. Mani-Levitska, C. Zong \cite{DLMZ02} for $n=4$. The case 
when the $2n$ congruent spherical caps are symmetric (see Lemma~\ref{2nevencovlemma} (i)) 
should be known, but we could not find any reference for the cases $n\geq 5$. 

\begin{lemma}
\label{2nevencovlemma}
Let $n\geq 2$, let $t\in(0,(2\cdot 4^{n-2}\sqrt{(n-1)!})^{-1})$,  and 
let $u_1,\ldots,u_n\in S^{n-1}$.
\begin{description}
\item{{\rm (i)}} If there exist $i<j$ such that $|\langle u_i,u_j\rangle|\geq \sin t$, 
then there exists  $u\in S^{n-1}$ such that
$$
|\langle u_i,u\rangle|\leq \frac1{\sqrt{n}}- \frac{t}{4n^{3/2}}
\mbox{ \ for $i=1,\ldots,n$}.
$$
\item{{\rm (ii)}} If $|\langle u_i,u_j\rangle|\leq \sin t$ for all $i<j$, then there 
exists a cross measure $\nu$ such that
$$
\delta_H({\rm supp}\,\nu,\{\pm u_1,\ldots,\pm u_n\})\leq 4^{n-2}\sqrt{(n-1)!}\cdot t.
$$
\end{description}
\end{lemma}
\proof For the proof of (i) we may assume that  $|\langle u_1,u_2\rangle|\geq \sin t$. 
We construct sequences $a_2,\ldots, a_n>0$
and $w_1,\ldots,w_n\in S^{n-1}$ such that $w_i\in{\rm lin}\{u_1,\ldots,u_i\}$,
and possibly after exchanging some of the vectors $u_i$ by $-u_i$, we have
$$
\langle w_i,u_j\rangle =a_i
\mbox{ \ for $i=1,\ldots,n$ and $j=1,\ldots,i$}.
$$
More precisely, let $w_1=u_1$, and if $i\in \{2,\ldots,n\}$ and $w_1,\ldots,w_{i-1}$ have 
already been determined, then we choose
the direction of $u_i$  in such a way that $\langle u_i,w_{i-1}\rangle\leq 0$. This 
algorithm determines $a_2,\ldots, a_n>0$
and $w_1,\ldots,w_n\in S^{n-1}$, and subsequently we prove that
\begin{equation}
\label{aiwi}
\langle w_i,u_j\rangle =a_i\leq \frac1{\sqrt{i}}-\frac{t}{4i^{3/2}}
\mbox{ \ for $i=2,\ldots,n$ and $j=1,\ldots,i$}.
\end{equation}

To verify (\ref{aiwi}), we use the elementary fact that if $o$ is a vertex 
of a triangle and if the 
two sides meeting at $o$ are of length $a$ and $b$ and enclose an angle $\gamma$, 
then the distance of $o$ from the line of the third side is 
\begin{equation}
\label{triangleheight}
h=\frac{ab\sin \gamma}{\sqrt{a^2+b^2-2ab\cos\gamma}}.
\end{equation}
In addition, we use that if $f(a)=\frac{a}{\sqrt{1+a^2}}$ for $a\in(0,s)$ and $s>0$, then
\begin{equation}
\label{aider}
f'(a)=\frac{1}{(1+a^2)^{3/2}}>\frac1{(1+s^2)^{3/2}}.
\end{equation}

We start with the case $i=2$. Since $\langle u_1,u_2\rangle \le 0$, we have
$\angle(u_1,u_2)\geq \frac{\pi}2+t$ and  $w_2=(u_1+u_2)/\|u_1+u_2\|$. Therefore
(\ref{triangleheight}) yields that
$$
\langle w_2,u_1\rangle=\langle w_2,u_2\rangle \le 
\frac{\cos t}{\sqrt{2+2\sin t}}<\frac1{\sqrt{2}}\cdot \frac1{\sqrt{1+\sin t}}<
\frac1{\sqrt{2}}\cdot \left(1-\frac{\sin t}4\right)<\frac1{\sqrt{2}}-\frac{t}{8\sqrt{2}}.
$$
Next assume that $2\leq i<n$ and (\ref{aiwi}) holds. We observe that 
$a_iw_i\in {\rm aff}\{u_1,\ldots,u_i\}$ and $a_{i+1}$ is 
the distance of $o$ from ${\rm aff}\{u_1,\ldots,u_{i+1}\}$, which is then 
at most the distance of $o$ from 
${\rm aff}\{a_iw_i,u_{i+1}\}$, that is in turn the height of the triangle 
$[o,a_iw_i,u_{i+1}]$ corresponding to $o$. Since
$\langle u_{i+1},w_i\rangle\leq 0$,  we deduce first from (\ref{triangleheight}),
then  from (\ref{aider}) with $a_i<s=\frac1{\sqrt{i}}$ that

$$
a_{i+1}\leq \frac{a_i}{\sqrt{1+a_i^2}}=f(a_i)< f(s)-\frac{t}{4i^{3/2}(1+s^2)^{3/2}}=
\frac1{\sqrt{i+1}}-\frac{t}{4(i+1)^{3/2}}.
$$
Finally, (\ref{aiwi}) yields  (i) with $u=w_n$. 

For (ii), let $v_1,\ldots,v_n$ be an orthonormal basis of $\R^n$ such that
$v_i\in{\rm lin}\{u_1,\ldots,u_i\}$ and $\langle v_j,u_j\rangle\geq 0$ for
$j=1,\ldots,n$, and hence $v_1=u_1$. We verify that
\begin{equation}
\label{goodbase}
\angle(v_i,u_i)\leq 4^{i-2}\sqrt{(i-1)!}\cdot t
\mbox{ \ for $i=2,\ldots,n$}
\end{equation}
by induction on $i=2,\ldots,n$.

If $i=2$, then readily $\angle(v_2,u_2)\leq t$. If (\ref{goodbase}) holds for all $j\le i$ for some  $i\in\{2,\ldots,n-1\}$, then 
$$
\left|\angle(u_{i+1},v_j)-\frac{\pi}2\right|\leq 
\left|\angle(u_{i+1},u_j)-\frac{\pi}2\right|+\angle(u_j,v_j)<2\cdot 4^{i-2}\sqrt{(i-1)!}\cdot t
$$
for $j=1,\ldots,i$.
In other words, $\langle u_{i+1},v_j\rangle<2\cdot 4^{i-2}\sqrt{(i-1)!}\cdot t$ 
for $j=1,\ldots,i$, which in turn yields that 
$$
\sin\angle(u_{i+1},v_{i+1})=\sqrt{\sum_{j=1}^i \langle u_{i+1},v_j\rangle^2}\le 2\cdot 4^{i-2}\sqrt{(i-1)!}\sqrt{i}\cdot t
=2\cdot 4^{i-2}\sqrt{i!}\cdot t.
$$
Thus we conclude $\angle(u_{i+1},v_{i+1})<4^{i-1}\sqrt{i!}\cdot t$.\proofbox

Lemma~\ref{2nevencovlemma} yields the following statement with factor
$4n^{3/2}\cdot 4^{n-2}\sqrt{(n-1)!}<4^nn!$.

\begin{coro}
\label{2nevencovcoro}
Let $n\geq 2$, let $t\in(0,\frac1{4^nn!})$,  and let $u_1,\ldots,u_n\in S^{n-1}$.
If 
$$
\Omega\left(u,\arccos\left(\frac1{\sqrt{n}}-t\right)\right)\cap\{\pm u_1,\ldots,\pm u_n\}\neq\emptyset
$$
for any $u\in S^{n-1}$, then there exists a cross measure $\nu$ such that
$$
\delta_H({\rm supp}\,\nu,\{\pm u_1,\ldots,\pm u_n\})\leq 4^nn!\cdot t.
$$
\end{coro}
{\bf Remark} The condition in Corollary \ref{2nevencovcoro} is equivalent to saying that 
$\Omega\left(\pm u_i,\arccos\left(\frac1{\sqrt{n}}-t\right)\right)$, $i=1,\ldots,n$, cover $S^{n-1}$.

\section{The volume of $Z^*_p$}
\label{secsymmetric}

In this section, we prove the stability result for the volume of $Z^*_p$, 
which is stated in Theorem~\ref{Zpmustab}. The remaining 
part of this theorem is established in Section~\ref{secLpzonoids}. 

The main ingredient for the proof in this section 
is stated as Proposition~\ref{volZ*pstab}.   
We start with preparatory claim.

\begin{claim}
\label{volumeuu_0}
For $u,u_0\in S^{n-1}$ with $\langle u,u_0\rangle\geq 0$, we have $V(\Xi_{u,u_0})\geq \kappa_n/240^n$,
where
\begin{equation}\label{volboundref}
\Xi_{u,u_0}=\left\{y\in 0.1\,B^n:\, \langle y,u\rangle\geq \frac1{30}, \, 
\langle y,u_0\rangle\geq \frac1{30}, \, \langle y,u-u_0\rangle\geq \frac{\|u-u_0\|}{120}
\right\}.
\end{equation}
\end{claim}
\proof Let $\gamma$ be  half of the angle of $u$ and $u_0$, and hence $\gamma\in [0,\frac{\pi}4]$.
The set
$$
\Xi_0=\left\{y\in  0.1\,B^n: \langle y,u\rangle\geq \frac1{30},\;\langle y,u_0\rangle\geq \frac1{30}\right\}
$$
contains a ball of radius $r$ with center $\frac{0.1-r}{\|u+u_0\|}\,(u+u_0)$ provided that
$$
(0.1-r)\cos\gamma\geq \mbox{$\frac1{30}$}+r.
$$
Since $\cos\gamma\geq 1/\sqrt{2}$, we may choose
$$
r=\frac{0.1-(\sqrt{2}/30)}{\sqrt{2}+1}>\frac1{60}.
$$
Therefore $\Xi_{u,u_0}$ contains a ball of radius $r/4>1/240$.
\proofbox

\begin{prop}
\label{volZ*pstab}
If $p\in[1,\infty)\setminus \{2\}$,  $\mu$ is an even discrete isotropic measure on $S^{n-1}$, and
$$
V(Z^*_p(\mu))\ge (1-\varepsilon) V(Z^*_p(\nu_n))
$$
for some $\varepsilon\in(0,1)$, then there exists
a cross measure $\nu$ on $S^{n-1}$  such that 
$$
\delta_W(\mu,\nu)\leq  n^{cn^3}\max\{|p-2|^{-\frac{2}3},1\}\cdot \varepsilon^{\frac{1}{3}}
$$
for some absolute constant $c>0$.
\end{prop}
\proof  What we actually prove is that for any $0<\eta<\beta^n/(2n)$, we have  
\begin{equation}
\label{volZ*pstab1}
V(Z^*_p(\mu))< (1-n^{-cn^3}\min\{(p-2)^2,1\}\cdot\eta^3) V(Z^*_p(\nu_n))
\end{equation}
or there exists a cross measure $\nu$ satisfying
\begin{equation}
\label{volZ*pstab2}
\delta_W(\mu,\nu)\leq  n^{cn}\eta
\end{equation}
for some absolute constant $c>0$.

Let ${\rm supp}\,\mu=\{\bar{u}_1,\ldots,\bar{u}_{\bar{k}}\}$, 
and  let $\bar{c}_i=\mu(\{\bar{u}_i\})$. For $c_0=\min\{\bar{c}_i:\,i=1,\ldots,\bar{k}\}$
and $i=1,\ldots,\bar{k}$, we define $\bar{m}_i=\min\{m\in\Z :\, m\geq 1\mbox{ and }\bar{c}_i/m\leq c_0\}$, and let
$k=\sum_{i=1}^{\bar{k}}\bar{m}_i$. We consider  $\xi:\{1,\ldots,k\}\to\{1,\ldots,\bar{k}\}$ such that
$\#\xi^{-1}(\{i\})=\bar{m}_i$ for $i=1,\ldots,\bar{k}$, and define
$$
u_i=\bar{u}_{\xi(i)}\mbox{ \ and \ }c_i=\bar{c}_{\xi(i)}/\bar{m}_{\xi(i)}
$$
for $i=1,\ldots,k$. The system $(u_1,\ldots,u_k,c_1,\ldots,c_k)$ is even (i.e.~origin symmetric) 
in the following sense: Any $u\in S^{n-1}$ occurs as $u_i$ exactly as many times as $-u$, and if $u_i=-u_j$, then $c_i=c_j$.

In particular,
$\sum_{i=1}^kc_i u_i\otimes u_i={\rm Id}_n$ and $\sum_{i=1}^kc_i=n$, 
and for any Borel $X\subseteq S^{n-1}$,  we have
$$
\mu(X)=\sum_{u_i\in X}c_i.
$$
The reason for the renormalization is that 
\begin{equation}
\label{cicj12}
c_0/2<c_i\leq c_0\mbox{ \ for $i=1,\ldots,k$}.
\end{equation}
In addition, let
$\varphi=\varphi_p$ be defined as in (\ref{phip}), let $g(t)=e^{-\pi t^2}$, and 
let $f_i=\varrho_p$, for $i=1,\ldots,k$.

We define the map $\Theta:\,\R^n\to\R^n$ by
$$
\Theta(y)=\sum_{i=1}^k  c_i\varphi(\langle y,u_i\rangle)\, u_i,
$$
and hence the differential of $\Theta$ is
$$
d\Theta(y)=\sum_{i=1}^k c_i \varphi'(\langle y,u_i\rangle)
\, u_i\otimes u_i ,
$$
where $d\Theta(y)$ is positive definite, and
 $\Theta:\R^n\to\R^n$ is injective. 
Applying first  (\ref{volZ*p0}) and then (\ref{BLstep1}), we get
\begin{align}
\nonumber
V(Z_p^*(\mu))&\leq \frac{2^n\Gamma(1+\frac{1}p)^n}{\Gamma(1+\frac{n}p)}
\int_{\R^n}\left(\prod_{i=1}^kg(\varphi(\langle u_i,x\rangle))^{c_i}\right)
\left(\prod_{i=1}^k\varphi'(\langle u_i,x\rangle)^{c_i}\right)\,dx\\
\label{volZp*form}
&= \frac{2^n\Gamma(1+\frac{1}p)^n}{\Gamma(1+\frac{n}p)}
\int_{\R^n}\exp\left(-\pi\sum_{i=1}^k c_i\varphi(\langle u_i,x\rangle)^2\right)
\left(\prod_{i=1}^k\varphi'(\langle u_i,x\rangle)^{c_i}\right)\,dx.
\end{align}

For each fixed $y\in \R^n$, we estimate the product of the two
terms in (\ref{volZp*form}) after the integral sign. To estimate the first term  in (\ref{volZp*form}), we apply (\ref{RBLfunc})  with $\theta_i=\varphi(\langle y,u_i\rangle)$, $i=1,\ldots,k$, and hence the definition
of $\Theta$ yields
\begin{equation}
\label{firsttermZp*}
\exp\left(-\pi \sum_{i=1}^k c_i\varphi(\langle y,u_i\rangle)^2  \right)\leq
 e^{-\pi\|\Theta(y)\|^2}.
\end{equation}

To estimate the second term, we apply Lemma~\ref{Ball-Barthe-stab}
with  $v_i=\sqrt{c_i}\cdot u_i$
and $t_i=\varphi'(\langle y,u_i\rangle)$, at each $y\in\R^n$, and
write $\theta^*(y)$ and $t_0(y)$ to denote
the corresponding $\theta^*\geq 1$ and $t_0>0$. In particular, 
if $ \{i_1,\ldots,i_n\}\subseteq\{1,\ldots,k\}$ and $y\in\R^n$, then we set
\begin{equation}\label{labelA}
\NN(i_1,\ldots,i_n;y)=c_{i_1}\cdots c_{i_n} 
\det[u_{i_1},\ldots,u_{i_n}]^2
\left(\frac{\sqrt{\varphi'(\langle y,u_{i_1}\rangle)\cdots 
\varphi'(\langle y,u_{i_n}\rangle)}}{t_0(y)}-1\right)^2.
\end{equation}
Therefore, for
\begin{equation}
\label{thetayZp*}
\theta^*(y)=  1+\frac12\sum_{1\leq i_1<\ldots<i_n\leq k}
\NN(i_1,\ldots,i_n;y)
\end{equation}
Lemma~\ref{Ball-Barthe-stab} yields that
\begin{equation}
\label{secondtermZp*}
\prod_{i=1}^k\varphi'(\langle y,u_i\rangle)^{c_i} \leq
\theta^*(y)^{-1} \det\left( d\Theta(y)\right).
\end{equation}
From (\ref{firsttermZp*}) and (\ref{secondtermZp*}), we conclude that 
\begin{equation}
\label{volZ*pstabstep1}
 V(Z_p^*(\mu))\leq \frac{2^n\Gamma(1+\frac{1}p)^n}{\Gamma(1+\frac{n}p)}
\int_{\R^n}\theta^*(y)^{-1} e^{-\pi\|\Theta(y)\|^2}\det \left(d\Theta(y)\right)\, dy.
\end{equation}

To provide a lower bound for $\theta^*(y)$, we use  (\ref{phipsiderest}) and \eqref{leqone}, hence
\begin{equation}
\label{phipsiderest-phip}
\frac1{3.1}<\varphi'(s)<3.1\mbox{ and } \varphi(s)<1
\mbox{ \ for $p\in[1,\infty]$ and $s\in[0,\frac1{3.1}]$.}
\end{equation}
We consider the vectors $v_1,\ldots,v_n\in S^{n-1}$ provided by Lemma~\ref{nsmallballs} 
such that 
\begin{eqnarray}
\label{Omegamu}
\mu(\Omega(v_i,\beta))&>&\beta^n\mbox{ \ \  for $i=1,\ldots,n$};\nonumber\\
\label{Omegadet}
|\det[w_1,\ldots,w_n]|&\geq& 2n\beta
\mbox{ \ \  for 
$w_i\in \Omega(v_i,\beta)$ and $i\in\{1,\ldots,n\}$};\label{reff83}\\
\label{betavalue}
\beta&=&2^{-(n+1)}n^{-(n+1)/2}.\nonumber
\end{eqnarray}
 The remaining discussion is split into three cases, where the first two correspond to the two cases in
Lemma~\ref{ntinyballs}.  \\

\noindent{\bf Case 1}\; {\em There exist $l\in\{1,\ldots,n\}$ and $\Psi_1,\Psi_2\subseteq \Omega(v_l,\beta)$ such that
\begin{eqnarray*}
\mu(\Psi_j)&\geq&\frac{\beta^n}{4n} \mbox{ \ for $j=1,2$},\text{ and}\\
\|a_1-a_2\|&\geq &\frac{\eta}{\sqrt{n}}
\mbox{ \ for $a_1\in\Psi_1$ and $a_2\in \Psi_2$}.
\end{eqnarray*} }
In this case, we prove
\begin{equation}
\label{Zp*case1}
 V(Z_p^*(\mu))<
\frac{2^n\Gamma(1+\frac{1}p)^n}{\Gamma(1+\frac{n}p)}
(1-n^{-cn^3}\min\{(p-2)^2,1\}\cdot\eta^2)
\end{equation}
for some absolute constant $c>0$. 

We may assume that $l=n$. For $j=1,2$, let 
$$
\Pi_j=\{i\in\{1,\ldots,k\}:\,u_i\in \Psi_j\}\neq \emptyset.
$$
 Possibly after interchanging the roles of $\Psi_1$ and $\Psi_2$, we may assume that
$\#\Pi_1\leq \#\Pi_2$. Let
$$
\tau:\Pi_1\to\Pi_2 \mbox{ \ be an injective map}.
$$
Given $u_{i_j}\in \Omega(v_j,\beta)$ for $j=1,\ldots,n-1$ and $u_{i_n}\in \Psi_1$, we have
have $u_{\tau(i_n)}\in \Psi_2$, 
and (\ref{cicj12}) and (\ref{Omegadet}) yield
\begin{equation}
\label{prodlower}
\left.
\begin{array}{r}
c_{i_1}\cdots c_{i_{n-1}}\cdot c_{i_n}\det[u_{i_1},\ldots,u_{i_n}]^2\\[1ex]
c_{i_1}\cdots c_{i_{n-1}}\cdot c_{\tau(i_n)}\det[u_{i_1},\ldots,u_{i_{n-1}},u_{\tau(i_n})]^2
\end{array}
\right\}\geq 4n^2\beta^2
c_{i_1}\cdots c_{i_{n-1}}\cdot(c_{i_n}/2).
\end{equation}
Since $\beta<\pi/4$, we have $\langle u_{i_n},u_{\tau(i_n)}\rangle >0$ if $u_{i_n}\in \Psi_1$. 
Claim~\ref{volumeuu_0} shows that $V(\Xi_{u,u_0})\geq \kappa_n/240^n$ 
for $u,u_0\in S^{n-1}$ with $\langle u,u_0\rangle\geq 0$, 
where $\Xi_{u,u_0}$ is defined in \eqref{volboundref}. 
In particular, if $y\in \Xi_{u_{i_n},u_{\tau(i_n)}}$, then
\begin{eqnarray*}
\langle y,u_{i_n}\rangle,\langle y,u_{\tau(i_n)}\rangle&<&\frac18,\text{ and}\\
\langle y,u_{i_n}\rangle-\langle y,u_{\tau(i_n)}\rangle&=&
\langle y,u_{i_n}-u_{\tau(i_n)}\rangle\geq \frac{\eta}{120\sqrt{n}}.
\end{eqnarray*}
Next, $\varphi''$ is continuous, and  Lemma~\ref{trasportation-map-phi} implies that if
$t\in[\frac1{30},0.1]$, then
\begin{equation}
\label{phidoubleder}
|\varphi''(t)|\geq\left\{
\begin{array}{lcll}
\frac{|p-2|}{48}\left(\frac1{30}\right)^{1.3}&>&\frac{|p-2|}{2^{12}}&
\mbox{ \ if $p\in[1,3]\setminus\{2\}$},\\[1ex]
0.2\left(\frac1{30}\right)^{1.3}&>&2^{-9}&
\mbox{ \ if $p>3$}.
\end{array} \right.
\end{equation}
Therefore,
$$
|\varphi'(\langle y,u_{i_n}\rangle)-\varphi'(\langle y,u_{\tau(i_n)}\rangle)|\ge 
\left\{
\begin{array}{lll}
\frac{|p-2|}{2^{12}120\sqrt{n}}\,\eta > &
 \frac{|p-2|}{2^{19}\sqrt{n}}\,\eta&\mbox{ \ if $p\in[1,3]\setminus \{2\}$},\\[1ex]
\frac1{2^{9}120\sqrt{n}}\,\eta>& \frac1{2^{19}\sqrt{n}}\,\eta&
\mbox{ \ if $p>3$}.
\end{array} \right.
$$
It follows from Lemma~\ref{xab} and $0< \varphi'(t)\leq 3.1$ for  $p\in[1,\infty)\setminus\{ 2\}$ and $t\in(0,0.1]$
(cf. (\ref{phipsiderest-phip})) that
\begin{align*}
&\left(\frac{\sqrt{\varphi'(\langle y,u_{i_1}\rangle)\cdots 
\varphi'(\langle y,u_{i_{n-1}}\rangle) \cdot
\varphi'(\langle y,u_{i_n}\rangle)}}{t_0(y)}-1\right)^2\\
&\quad +
\left(\frac{\sqrt{\varphi'(\langle y,u_{i_1}\rangle)\cdots  
\varphi'(\langle y,u_{i_{n-1}}\rangle) \cdot
\varphi'(\langle y,u_{\tau(i_n)}\rangle)}}{t_0(y)}-1\right)^2&\\
&\qquad\qquad\geq \frac{(\varphi'(\langle y,u_{i_n}\rangle)-\varphi'(\langle y,u_{\tau(i_n)}\rangle))^2}
{2(\varphi'(\langle y,u_{i_n}\rangle)+\varphi'(\langle y,u_{\tau(i_n)}\rangle))^2}
\ge \frac{\min\{1,(p-2)^2\}}{2^{45}n}\,\eta^2.
\end{align*}
Combining this estimate with \eqref{labelA} and (\ref{prodlower})  
implies that if   $p\in[1,\infty)\setminus \{ 2\}$ and $u_{i_j}\in \Omega(v_j,\beta)$ for $j=1,\ldots,n-1$, 
$u_{i_n}\in \Psi_1$ and $y\in \Xi_{u_{i_n},u_{\tau(i_n)}}$, then
\begin{align*}
&\NN(i_1,\ldots,i_{n-1},i_n;y)+\NN(i_1,\ldots,i_{n-1},\tau(i_n);y)\\
&\qquad\ge 4n^2\beta^2
c_{i_1}\cdots c_{i_{n-1}}\cdot(c_{i_n}/2)
\frac{\min\{1,(p-2)^2\}}{2^{45}n}\,\eta^2.
\end{align*}

If $u_{i_n}\in \Psi_1$ and $y\in\R^n$, then we define
$$
\varrho(i_n;y)=
\left\{
\begin{array}{ll}
0&\mbox{ \ if $y\not\in \Xi_{i_n,\tau(i_n)}$};\\
\frac{\beta^2n(p-2)^2}{2^{44}}\,\eta^2&\mbox{ \ if $y\in \Xi_{i_n,\tau(i_n)}$ and $p\in[1,3]\setminus\{2\}$};\\[1ex]
 \frac{\beta^2n}{2^{44}}\,\eta^2&
\mbox{ \ if $y\in \Xi_{i_n,\tau(i_n)}$ and $p>3$}.
\end{array} \right.
$$
In particular, if $u_{i_j}\in \Omega(v_j,\beta)$ for $j=1,\ldots,n-1$, 
$u_{i_n}\in \Psi_1$ and $y\in \R^n$, then
\begin{equation}
\label{alephest}
\NN(i_1,\ldots,i_{n-1},i_n;y)+\NN(i_1,\ldots,i_{n-1},\tau(i_n),y)\geq 
c_{i_1}\cdots c_{i_n}\varrho(i_n;y).
\end{equation}
Substituting (\ref{alephest}) into (\ref{thetayZp*}), and then using
(\ref{Omegamu}), 
 we see that if $y\in\R^n$, then
\begin{align*}
\theta^*(y)&\geq   1+\frac12\sum_{u_{i_j}\in \Omega(v_j,\beta),\;
j=1,\ldots,n-1\atop u_{i_n}\in \Psi_1}
c_{i_1}\cdots c_{i_{n-1}}\cdot c_{i_n}\varrho(i_n;y)\\
&= 1+\frac12\left(\prod_{j=1}^{n-1}\mu(\Omega(v_j,\beta))\right)
\sum_{u_{i_n}\in \Psi_1} c_{i_n}\varrho(i_n;y)\\
&\geq
1+\frac{\beta^{n(n-1)}}{2}\sum_{u_{i_n}\in \Psi_1} c_{i_n}\varrho(i_n;y).
\end{align*}
Here 
$$
\frac{\beta^{n(n-1)}}2\sum_{u_{i_n}\in \Psi_1} c_{i_n}\varrho(i_n;y)\leq
\frac{\beta^{n(n-1)}}2\,\mu(\Psi_1)\cdot\frac{\beta^2n}{2^{44}}\,\eta^2<1,
$$
and hence  if $y\in\R^n$, then
\begin{equation}
\label{theta*recipp}
\theta^*(y)^{-1}\leq 1-\frac{\beta^{n(n-1)}}4\sum_{u_{i_n}\in \Psi_1} c_{i_n}\varrho(i_n;y).
\end{equation}

We deduce from (\ref{volZ*pstabstep1}) and (\ref{theta*recipp})  that
\begin{align*}
 V(Z_p^*(\mu))&\leq \frac{2^n\Gamma(1+\frac{1}p)^n}{\Gamma(1+\frac{n}p)}
\int_{\R^n}\left(1-\frac{\beta^{n(n-1)}}4\sum_{u_{i_n}\in \Psi_1} c_{i_n}
\varrho(i_n;y)\right) e^{-\pi\|\Theta(y)\|^2}\det \left(d\Theta(y)\right)\,dy\\
&= \frac{2^n\Gamma(1+\frac{1}p)^n}{\Gamma(1+\frac{n}p)}
\int_{\R^n}e^{-\pi\|\Theta(y)\|^2}\det\left( d\Theta(y)\right)\,dy\\
&\qquad- \frac{2^n\Gamma(1+\frac{1}p)^n}{\Gamma(1+\frac{n}p)}\cdot
\frac{\beta^{n(n-1)}}4 \sum_{u_{i_n}\in \Psi_1} c_{i_n}
\int_{\R^n}\varrho(i_n;y)e^{-\pi\|\Theta(y)\|^2}\det\left( d\Theta(y)\right)\,dy.
\end{align*}
Here, we use that 
\begin{equation}\label{leqthanone}
\int_{\R^n}e^{-\pi\|\Theta(y)\|^2}\det\left( d\Theta(y)\right)\,dy \le 
\int_{\R^n} e^{-\pi\|z\|^2}\,dz=1.
\end{equation}
If $y\in \Xi_{i_n,\tau(i_n)}$, then 
 (\ref{firsttermZp*}), (\ref{secondtermZp*})  and (\ref{phipsiderest-phip}) yield that
\begin{eqnarray}
\label{firsttermZp*cons}
 e^{-\pi \|\Theta(y)\|^2} &\geq & 
\exp\left(-\pi \sum_{i=1}^k c_i\varphi(\langle y,u_i\rangle)^2  \right)>
\exp\left(-\pi \sum_{i=1}^k c_i \right)=e^{-\pi n},\\
\label{secondtermZp*cons}
\det\left( d\Theta(y)\right)&\geq &\prod_{i=1}^k\varphi'(\langle y,u_i\rangle)^{c_i}\geq
\prod_{i=1}^k3.1^{-c_i}=3.1^{-n}.
\end{eqnarray}
Therefore
$$
 V(Z_p^*(\mu))\leq \frac{2^n\Gamma(1+\frac{1}p)^n}{\Gamma(1+\frac{n}p)}
\Bigg(1-\sum_{u_{i_n}\in \Psi_1} c_{i_n}\frac{\beta^{n(n-1)}}4\cdot 
\frac{V(\Xi_{i_n,\tau(i_n)})}{(3.1e^\pi)^n}
\cdot \frac{\beta^2n\min\{(p-2)^2,1\}}{2^{44}}\cdot\eta^2 \Bigg).
$$
Since $V(\Xi_{i_n,\tau(i_n)})\geq \kappa_n/240^n$ if
$u_{i_n}\in \Psi_1$, according to  Claim~\ref{volumeuu_0}, and 
$$\sum_{u_{i_n}\in \Psi_1} c_{i_n}=\mu(\Psi_1)>\frac{\beta^n}{4n},$$ 
we conclude (\ref{Zp*case1}).\\

\noindent{\bf Case 2}\;  {\em There exists $q_i\in\Omega(v_i,\beta)$, for $i=1,\ldots,n$, such that
\begin{eqnarray}
&\mu(\Omega(q_i,\eta))\geq \frac{\beta^n}{4n}\text{ \ \rm for $i=1,\ldots,n$},\text{ \rm and}\label{case2a}\\
&\mu\left(\bigcup_{i=1}^n(\Omega(q_i,2\eta)\cup \Omega(-q_i,2\eta))\right)\le n-\eta.\label{case2b}
\end{eqnarray}}
In this case, we prove
\begin{equation}
\label{Zp*case2}
 V(Z_p^*(\mu))\le 
\frac{2^n\Gamma(1+\frac{1}p)^n}{\Gamma(1+\frac{n}p)}
(1-n^{-cn^3}\min\{(p-2)^2,1\}\cdot\eta^3)
\end{equation}
for some absolute constant $c>0$. The argument is very similar to the one in Case 1.

Let
$$
\widetilde{\Psi}=S^{n-1}\setminus \left(\bigcup_{i=1}^n (\Omega(q_i,2\eta)\cup \Omega(-q_i,2\eta))\right).
$$
It follows from \eqref{reff83} that any $x\in\R^n$ can be written in the form
$$
x=\sum_{i=1}^n\lambda_i(x)q_i.
$$
Since $\mu(\widetilde{\Psi})\ge \eta$ by \eqref{case2b}, the triangle inequality ensures that there exists some $i\in\{1,\ldots,n\}$ satisfying $|\lambda_i(x)|\geq 1/n$. Thus we may reindex $q_1,\ldots,q_n$ in such a way that
\begin{equation}
\label{Psibig}
\mu(\Psi)\ge\frac{\eta}n \quad\mbox{ \ for }
\Psi=\{x\in \widetilde{\Psi}:\,|\lambda_n(x)|\geq 1/n \}.
\end{equation}
We deduce from (\ref{Omegadet}) that if $x\in \Psi$, then
$$
|\det[q_1,\ldots,q_{n-1},x]|\geq |\det[q_1,\ldots,q_{n-1},q_n]|/n\geq
 2\beta.
$$
Next, for $u_{i_j}\in \Omega(q_j,\eta)$ for $j=1,\ldots,n-1$, we apply Claim~\ref{largedeterminant} 
with $b_l=q_l$, $s_l=u_{i_l}-q_l$, for $l=1,\ldots,n-1$, 
 $b_n=x\in\Psi$, and $s_n=0$, where  
$$
|s_i|\le \eta\le \frac{\beta}{2n}=\frac{2\beta}{4n}\le \frac{1}{4n}|\det[q_1,\ldots,q_{n-1},x]|,\quad i=1,\ldots,n.
$$
Hence, 
\begin{equation}
\label{Psidet}
|\det[u_{i_1},\ldots,u_{i_{n-1}},x]|\geq \frac{1}{2} |\det[q_1,\ldots,q_{n-1},x]|\geq\beta.
\end{equation}

We observe that $\Psi=-\Psi$.
Thus, for  
$$
\Pi_2=\{i\in\{1,\ldots,k\}:u_i\in\Psi\},
$$
there exists $\Pi'\subseteq\Pi_2$ with $\#\Pi'=\frac12 \#\Pi_2$, and  a bijection
$\tilde{\tau}:\Pi'\to\Pi_2\setminus \Pi'$ such that if $i\in\Pi'$ then 
$u_{\tilde{\tau}(i)}=-u_i$.

Since $\eta<\beta^n$, \eqref{case2a} implies that  
$$
\sum_{u_i\in \Omega(q_n,\eta)} c_i=\mu(\Omega(q_n,\eta))\ge \frac{\beta^n}{4n}\ge\frac{\eta}{8n}.
$$ 
Thus we can find a minimal (with respect to inclusion) set $\Pi_1\subseteq\{1,\ldots,k\}$ such that 
$u_i\in\Omega(q_n,\eta)$ for $i\in\Pi_1$ and
\begin{equation}\label{Psi1big}
\sum_{i\in \Pi_1} c_i\ge\frac{\eta}{8n},
\end{equation}
By minimality and \eqref{cicj12} it follows that
$$
\frac{c_0}{2}\left(\#\Pi_1-1\right)\leq \frac{\eta}{8n}.
$$
Moreover, by  \eqref{Psibig} and again by \eqref{cicj12}, we have
$$
c_0\#\Pi_2\ge \sum_{j\in \Pi_2} c_j\ge \frac{\eta}n,
$$
and hence
$$
\frac{c_0}{8}\#\Pi_2\ge \frac{c_0}{2}\left(\#\Pi_1-1\right),
$$
which yields $\#\Pi_2\ge 4(\#\Pi_1-1)$ if $\#\Pi_1\ge 2$. In any case, we deduce 
that $\#\Pi_2\ge 2\#\Pi_1$.

We conclude that there exists
an injective map  $\tau:\Pi_1\to\Pi_2$ such that if $i\in\Pi_1$, then
\begin{equation}
\label{Psi1tauacute}
\langle u_i,u_{\tau(i)}\rangle \geq 0.
\end{equation}
In addition, if $i\in\Pi_1$, then $u_i\in\Omega(q_n,\eta)$ and 
 $u_{\tau(i)}\not\in\Omega(q_n,2\eta)$, and therefore
$$
\| u_i-u_{\tau(i)}\| \geq \frac{\eta}{2}.
$$

Given $u_{i_j}\in \Omega(q_j,\eta)$ for $j=1,\ldots,n-1$ and $i_n\in \Pi_1$, we have
have $\tau(i_n)\in \Pi_2$, 
and (\ref{cicj12}), (\ref{Omegadet}) and (\ref{Psidet}) yield
\begin{equation}
\label{prodlowerCase2}
\left.
\begin{array}{r}
c_{i_1}\cdots c_{i_{n-1}}\cdot c_{i_n}\det[u_{i_1},\ldots,u_{i_n}]^2\\[1ex]
c_{i_1}\cdots c_{i_{n-1}}\cdot c_{\tau(i_n)}\det[u_{i_1},\ldots,u_{i_{n-1}},u_{\tau(i_n})]^2
\end{array}
\right\}\geq \beta^2
c_{i_1}\cdots c_{i_{n-1}}\cdot(c_{i_n}/2).
\end{equation}
We deduce from (\ref{Psi1tauacute}) that  Claim~\ref{volumeuu_0} applies to
$\Xi_{u_{i_n},u_{\tau(i_n)}}$. In particular, we have 
$V(\Xi_{u_{i_n},u_{\tau(i_n)}})\geq \kappa_n/240^n$, and 
 if $y\in \Xi_{u_{i_n},u_{\tau(i_n)}}$, then
\begin{eqnarray*}
\langle y,u_{i_n}\rangle,\langle y,u_{\tau(i_n)}\rangle&<&\frac18;\\
\langle y,u_{i_n}\rangle-\langle y,u_{\tau(i_n)}\rangle&=&
\langle y,u_{i_n}-u_{\tau(i_n)}\rangle\geq \frac{\eta}{240}>\frac{\eta}{2^8}.
\end{eqnarray*}
It follows from (\ref{phidoubleder}) that
$$
|\varphi'(\langle y,u_{i_n}\rangle)-\varphi'(\langle y,u_{\tau(i_n)}\rangle)|\ge 
 \frac{\min\{|p-2|,1\}}{2^{20}}\cdot\eta.
$$
Since $0<\varphi'(t)\leq 3.1$ for $t\in(0,0.1]$, if $i_n\in\Pi_1$, then
$$
\frac{(\varphi'(\langle y,u_{i_n}\rangle)-\varphi'(\langle y,u_{\tau(i_n)}\rangle))^2}
{2(\varphi'(\langle y,u_{i_n}\rangle)+\varphi'(\langle y,u_{\tau(i_n)}\rangle))^2}\ge 
 \frac{\min\{(p-2)^2,1\}}{2^{47}}\cdot\eta^2.
$$
Thus combining Lemma~\ref{xab} and  (\ref{prodlowerCase2}), we obtain that
 if $u_{i_j}\in \Omega(v_j,\beta)$ for $j=1,\ldots,n-1$, 
$i_n\in \Pi_1$ and $y\in \Xi_{u_{i_n},u_{\tau(i_n)}}$, then
$$
\NN(i_1,\ldots,i_{n-1},i_n;y)+\NN(i_1,\ldots,i_{n-1},\tau(i_n);y)\geq 
\frac{\beta^2 c_{i_1}\cdots c_{i_n}}2\cdot 
\frac{\min\{(p-2)^2,1\}}{2^{47}}\cdot\eta^2.
$$
If $i_n\in \Pi_1$ and $y\in\R^n$, then we define
$$
\varrho(i_n;y)=
\left\{
\begin{array}{ll}
0&\mbox{ \ if $y\not\in \Xi_{i_n,\tau(i_n)}$}\\
\frac{\beta^2\min\{(p-2)^2,1\}}{2^{48}}\cdot\eta^2&\mbox{ \ if $y\in \Xi_{i_n,\tau(i_n)}$}.
\end{array} \right.
$$
In particular, if $u_{i_j}\in \Omega(v_j,\beta)$ for $j=1,\ldots,n-1$, 
$i_n\in \Pi_1$ and $y\in \R^n$, then
\begin{equation}
\label{alephestCase2}
\NN(i_1,\ldots,i_{n-1},i_n;y)+\NN(i_1,\ldots,i_{n-1},\tau(i_n),y)\geq 
c_{i_1}\cdots c_{i_n}\varrho(i_n;y).
\end{equation}
Substituting (\ref{alephestCase2}) into (\ref{thetayZp*}) and then using
(\ref{Omegamu}), 
 we obtain for $y\in\R^n$ that
\begin{align*}
\theta^*(y)&\geq   1+\frac12\sum_{u_{i_j}\in \Omega(v_j,\beta),\;j=1,\ldots,n-1\atop i_n\in \Pi_1}
c_{i_1}\cdots c_{i_{n-1}}\cdot c_{i_n}\varrho(i_n;y)\\
&=1+\frac12\left(\prod_{j=1}^{n-1}\mu(\Omega(v_j,\beta))\right)
\sum_{i_n\in \Pi_1} c_{i_n}\varrho(i_n;y)\\
&\geq
1+\frac{\beta^{n(n-1)}}2\sum_{i_n\in \Pi_1} c_{i_n}\varrho(i_n;y).
\end{align*}
Similarly as before, we have
$$
\frac{\beta^{n(n-1)}}2\sum_{i_n\in \Pi_1} c_{i_n}\varrho(i_n;y)\leq
\frac{\beta^{n(n-1)}}2\,\mu(\Psi_1)\cdot\frac{\beta^2n}{2^{48}}\cdot\eta^2<1,
$$
and hence 
\begin{equation}
\label{theta*recippCase2}
\theta^*(y)^{-1}\leq 1-\frac{\beta^{n(n-1)}}4\sum_{i_n\in \Pi_1} c_{i_n}\varrho(i_n;y).
\end{equation}

We deduce from (\ref{volZ*pstabstep1}) and (\ref{theta*recippCase2})  that 
\begin{align*}
 V(Z_p^*(\mu))&\leq  \frac{2^n\Gamma(1+\frac{1}p)^n}{\Gamma(1+\frac{n}p)}
\int_{\R^n}e^{-\pi\|\Theta(y)\|^2}\det\left( d\Theta(y)\right)\,dy\\
&\qquad- \frac{2^n\Gamma(1+\frac{1}p)^n}{\Gamma(1+\frac{n}p)}\cdot
\frac{\beta^{n(n-1)}}4 \sum_{i_n\in \Pi_1} c_{i_n}
\int_{\R^n}\varrho(i_n;y)e^{-\pi\|\Theta(y)\|^2}\det\left( d\Theta(y)\right)\,dy.
\end{align*}
Now we use again \eqref{leqthanone} as well as the estimates 
(\ref{firsttermZp*cons}) and (\ref{secondtermZp*cons})
if $y\in \Xi_{i_n,\tau(i_n)}$.
Therefore
$$
 V(Z_p^*(\mu))\leq \frac{2^n\Gamma(1+\frac{1}p)^n}{\Gamma(1+\frac{n}p)}
\left(1-\sum_{i_n\in \Pi_1} c_{i_n}\frac{\beta^{n(n-1)}}4\cdot 
\frac{V(\Xi_{i_n,\tau(i_n)})}{(3.1e^\pi)^n}
\cdot \frac{\beta^2\min\{(p-2)^2,1\}}{2^{48}}\cdot\eta^2 \right)
$$
Since $V(\Xi_{i_n,\tau(i_n)})\geq \kappa_n/240^n$ if
$i_n\in \Pi_1$ and by (\ref{Psi1big}), 
we conclude (\ref{Zp*case2}).\\

\noindent{\bf Case 3}\; {\em There exists $q_i\in\Omega(v_i,\beta)$, for $i=1,\ldots,n$, such that
$$
\mu\left(\bigcup_{i=1}^n(\Omega(q_i,2\eta)\cup \Omega(-q_i,2\eta))\right) > n-\eta.
$$}

In this case, we prove that there exists a cross measure $\nu$ such that
\begin{equation}
\label{Zp*case3}
\delta_W(\nu,\mu)\leq n^{cn}\eta
\end{equation}
for some absolute constant $c>0$.

We observe that $\frac12(1-n(\frac1{\sqrt{n}}-t)^2)>\eta$
for $t=2\eta$, since $\eta<1/(2n)$. Thus Claim~\ref{isotropiccap} yields that
$\Omega(u,\arccos(\frac1{\sqrt{n}}-2\eta))$ intersects
$\cup_{i=1}^n\Omega(\pm q_i,2\eta)$ for any $u\in S^{n-1}$. In turn, we deduce that
$$
\Omega\left(u,\arccos\left(\frac1{\sqrt{n}}-4\eta\right)\right)\cap\{\pm q_1,\ldots,\pm q_n\}\neq\emptyset
$$
for any $u\in S^{n-1}$, since $4\eta<1/(4^n n!)$. 
Therefore Corollary~\ref{2nevencovcoro} implies that there exists a cross measure $\nu$ such that
$$
\delta_H(\text{\rm supp }\nu,\{\pm q_1,\ldots,\pm q_n\})\leq 4^nn!\cdot 4\eta.
$$
In particular, (\ref{Zp*case3}) follows from Lemma~\ref{Hausdorffcrosserror}.

According to
Lemma~\ref{ntinyballs},  Cases 1, 2 and 3 cover all possible even isotropic measure $\mu$. 
Thus, we have proved (\ref{volZ*pstab1}) in Cases 1 and 2, and
(\ref{volZ*pstab2}) in Case 3.
\proofbox

\noindent {\bf Proof of Theorem~\ref{Zpmustab} in the case of $Z^*_p(\mu)$: } Let 
 $p\in[1,\infty)\setminus \{2\}$, and let $\mu$ be a discrete even isotropic measure on $S^{n-1}$. Assume 
that  $\delta_{\rm WO}(\mu,\nu_n)\geq\varepsilon>0$ for some $\varepsilon\in(0,1)$. Then
Proposition~\ref{volZ*pstab} yields that
\begin{equation}
\label{Z*pmustabdiscrete}
V(Z^*_p(\mu))\le (1-\gamma\varepsilon^3) V(Z^*_p(\nu_n)),
\end{equation}
where $\gamma=n^{-cn^3}\min\{|p-2|^2,1\}$ for an absolute constant $c>0$.
Since any even isotropic measure can be weakly approximated by discrete even  isotropic measures 
(see, for instance,  F. Barthe \cite{Bar04}), we conclude (\ref{Z*pmustabdiscrete}), and in turn
Theorem~\ref{Zpmustab} in the case of $Z^*_p(\mu)$, for any even isotropic measure $\mu$ on $S^{n-1}$ 
and $p\in[1,\infty)\setminus\{2\}$.

 Since for any isotropic measure $\mu$, we have
$$
\lim_{p\to\infty} Z^*_p(\mu)=Z^*_{\infty}(\mu),
$$
and the factor $\gamma$ in (\ref{Z*pmustabdiscrete}) is independent of  $p\in (2,\infty)$,
we deduce the case $p=\infty$ as well. \proofbox

\section{The case of the $L_p$ zonoids in Theorem~\ref{Zpmustab}}
\label{secLpzonoids}

The proof of Theorem~\ref{Zpmustab} for $V(Z_p(\mu))$ is analogous to the argument for $V(Z_p^*(\mu))$.
In particular, we may assume again that $\mu$ is a discrete even isotropic measure, and
$p\in(1,\infty)\setminus \{2\}$. Let $p^*\in(1,\infty)$ be defined by $\frac1p+\frac1{p^*}=1$.
We prove that if $\eta\in(0,1)$, then  
\begin{equation}
\label{volZp*stab1}
V(Z_{p^*}(\mu))> (1-n^{-cn^3}\min\{(p-2)^2,1\}\cdot\eta^3) V(Z_{p^*}(\nu_n)) 
\end{equation}
or there exists a cross measure $\nu$ satisfying
\begin{equation}
\label{volZp*stab2}
\delta_W(\mu,\nu)\leq  n^{cn}\eta
\end{equation}
for some absolute constant $c>0$. Since if $p\in[\frac32,3]$, then $p^*\in[\frac32,3]$ and
$|p-2|/2\leq |p^*-2|\leq 2|p-2|$, (\ref{volZp*stab1}) and (\ref{volZp*stab2}) yield Theorem~\ref{Zpmustab} for $V(Z_p(\mu))$.

Again, let ${\rm supp}\,\mu=\{\bar{u}_1,\ldots,\bar{u}_{\bar{k}}\}$, 
and  let $\bar{c}_i=\mu(\{\bar{u}_i\})$. For $c_0=\min\{\bar{c}_i:\,i=1,\ldots,\bar{k}\}$
and $i=1,\ldots,\bar{k}$, we define $\bar{m}_i=\min\{m\in\Z :\, m\geq 1\mbox{ and }\bar{c}_i/m\leq c_0\}$, and let
$k=\sum_{i=1}^{\bar{k}}\bar{m}_i$. We consider  $\xi:\{1,\ldots,k\}\to\{1,\ldots,\bar{k}\}$ such that
$\#\xi^{-1}(\{i\})=\bar{m}_i$ for $i=1,\ldots,\bar{k}$, and define
$$
u_i=\bar{u}_{\xi(i)}\mbox{ \ and \ }c_i=\bar{c}_{\xi(i)}/\bar{m}_{\xi(i)}
$$
for $i=1,\ldots,k$. 

In particular,
$\sum_{i=1}^kc_i u_i\otimes u_i={\rm Id}_n$ and $\sum_{i=1}^kc_i=n$, and for any Borel $X\subseteq S^{n-1}$,  we have
$$
\mu(X)=\sum_{u_i\in X}c_i.
$$
Again, we obtain
$$
c_0/2<c_i\leq c_0\mbox{ \ for $i=1,\ldots,k$}.
$$
In addition, let
$\psi=\psi_p$ be defined as in (\ref{psip}), let $g(t)=e^{-\pi t^2}$, and let 
$f_i=\varrho_p$, for $i=1,\ldots,k$. 

We define the map $\Psi:\,\R^n\to\R^n$ by
$$
\Psi(y)=\sum_{i=1}^k  c_i\psi(\langle y,u_i\rangle)\, u_i.
$$
Its differential 
$$
d\Psi(y)=\sum_{i=1}^k c_i \psi'(\langle y,u_i\rangle)\,  u_i\otimes u_i 
$$
is positive definite, and
 $\Psi:\R^n\to\R^n$ is injective.

It follows by first applying  (\ref{volZp*RBL}), and then (\ref{RBLstep1}), that
\begin{align*}
V(Z_{p^*}(\mu))&\geq V(M_p(\mu))=\frac{2^n\Gamma(1+\frac{1}p)^n}{\Gamma(1+\frac{n}p)}
\int_{\R^n}^\ast\sup_{x=\sum_{i=1}^kc_i\theta_iu_i}  \prod_{i=1}^kf_i(\theta_i)^{c_i}\,dx\\
&\geq  \frac{2^n\Gamma(1+\frac{1}p)^n}{\Gamma(1+\frac{n}p)}
\int_{\R^n}\left(\prod_{i=1}^kf_i(\psi(\langle u_i,y\rangle))^{c_i} \right)
\det\left(\sum_{i=1}^kc_i\psi'(\langle u_i,y\rangle )\,u_i\otimes u_i\right)\,dy.
\end{align*}

To estimate the second term, we apply Lemma~\ref{Ball-Barthe-stab}
with  $v_i=\sqrt{c_i}\cdot u_i$
and $t_i=\psi'(\langle y,u_i\rangle)$ at each $y\in\R^n$, and
write $\theta^*(y)$ and $t_0(y)$ to denote
the corresponding $\theta^*\geq 1$ and $t_0$. In particular, 
if $ \{i_1,\ldots,i_n\}\subseteq\{1,\ldots,k\}$ and $y\in\R^n$, then we now set
$$
\NN(i_1,\ldots,i_n;y)=c_{i_1}\cdots c_{i_n}
\det[u_{i_1},\ldots,u_{i_n}]^2
\left(\frac{\sqrt{\psi'(\langle y,u_{i_1}\rangle)\cdots 
\psi'(\langle y,u_{i_n}\rangle)}}{t_0(y)}-1\right)^2.
$$
Therefore, using again the notation
$$
\theta^*(y)=  1+\frac12\sum_{1\leq i_1<\ldots<i_n\leq k}\NN(i_1,\ldots,i_n;y),
$$
Lemma~\ref{Ball-Barthe-stab} and (\ref{masstransS})  lead to
\begin{align*}
V(Z_{p^*}(\mu))&\geq  \frac{2^n\Gamma(1+\frac{1}p)^n}{\Gamma(1+\frac{n}p)}
\int_{\R^n}\theta^*(y)\left(\prod_{i=1}^kf_i(\psi(\langle u_i,y\rangle))^{c_i} \right)
\left(\prod_{i=1}^k\psi'(\langle u_i,y\rangle )^{c_i}\right)\,dy\\
&=\frac{2^n\Gamma(1+\frac{1}p)^n}{\Gamma(1+\frac{n}p)}
\int_{\R^n}\theta^*(y) \left(\prod_{i=1}^k g(\langle u_i,y\rangle )^{c_i}\right)\,dy\\
&=\frac{2^n\Gamma(1+\frac{1}p)^n}{\Gamma(1+\frac{n}p)}
\int_{\R^n}\theta^*(y) e^{-\pi\|y\|^2}\,dy.
\end{align*}
Now (\ref{volZp*stab1}) and (\ref{volZp*stab2}), and hence 
Theorem~\ref{Zpmustab} for $V(Z_p(\mu))$, can be proved as
(\ref{volZ*pstab1}) and (\ref{volZ*pstab2}) in Proposition~\ref{volZ*pstab}
were proved following (\ref{secondtermZp*}).

\section{Stability of the  reverse isoperimetric inequality in the origin symmetric case}
\label{secreverseisostab}

In this section, we turn to the proofs of Corollary \ref{Zinftymustab} and 
of Theorems \ref{vol} and \ref{BM}. 

We may assume that the facets of the cube $W^n$ touch $B^n$ in the support of the 
reference  cross measure $\nu_n$, where 
${\rm supp}\,\nu_n=\{\pm e_1,\ldots,\pm e_n\}$. 

\begin{lemma}
\label{cubesandwich}
If $\mu$ is an even measure on $S^{n-1}$ such that
$\delta_H({\rm supp}\,\mu,{\rm supp}\,\nu_n)<\alpha$ for some 
$\alpha\in(0,\frac1{3n})$, then $e^{-n\alpha}W^n
\subseteq Z^*_\infty(\mu)\subseteq e^{2n\alpha}W^n$.
\end{lemma}
\proof 
First, we show that $Z^*_\infty(\mu)\subseteq e^{2n\alpha}W^n$. 
For this, let $x\in\R^n\setminus e^{2n\alpha}W^n$. 
Clearly, we may assume that $x_1=\max\{|x_1|,\ldots,|x_n|\}$. 
It follows that there is some $i\in\{1,\ldots,n\}$ such that
\begin{equation}\label{instar1}
x_1\ge |x_i|=|\langle x,e_i\rangle|>e^{2n\alpha}\ge 
\left(1-\frac{1}{2}\alpha^2-\sqrt{n-1}\,\alpha\right)^{-1}
\ge \left(\cos\alpha-\sqrt{n-1}\sin\alpha\right)^{-1},
\end{equation}
where we used that $\alpha\in(0,\frac1{3n})$ for the third inequality. Since 
$\delta_H({\rm supp}\,\mu,{\rm supp}\,\nu_n)<\alpha$, 
there is some $v\in {\rm supp}\,\mu$ such that 
$\angle (e_1,v)< \alpha$, hence
\begin{equation}\label{instar2}
\langle e_1,v\rangle > \cos\alpha,\qquad  
\sum_{i=2}^n |\langle e_i,v\rangle| < \sqrt{n-1}\sin\alpha.
\end{equation}
From \eqref{instar1} we deduce that
$$
\langle x,v\rangle\geq x_1\langle e_1,v\rangle-x_1\sum_{i=2}^n 
|\langle e_i,v\rangle|> x_1\left(\cos\alpha - \sqrt{n-1}\sin\alpha\right)>
1,
$$
and hence $x\notin Z^*_\infty(\mu)$. 

In order to show that $e^{-n\alpha}W^n\subseteq Z^*_\infty(\mu)$, we put 
$\varrho=(1+\sqrt{n-1}\sin\alpha)^{-1}$. Since $\varrho\ge (1+n\alpha)^{-1}\ge e^{-n\alpha}$, 
we have $e^{-n\alpha} W^n\subseteq \varrho W^n$, and it is sufficient to show that 
$\varrho W^n\subseteq  Z^*_\infty(\mu)$. For this, let $x\in \varrho W^n$, and let $v\in {\rm supp}\,\mu$ 
be arbitrary. Then there is some $i\in\{1,\ldots,n\}$ such that $\angle (e_i,v)< \alpha$ or 
$\angle (-e_i,v)< \alpha$. We may assume that $i=1$. Hence \eqref{instar2} is available again. 
Then $x=x_1e_1+\ldots+x_ne_n$ with $|x_i|\le \varrho$ satisfies
$$
\langle x,v\rangle\leq \varrho\cdot 1+\varrho\sqrt{n-1}\sin\alpha = 1,
$$ 
which shows that $x\in Z^*_\infty(\mu)$.
\proofbox

For the proof of  Theorem~\ref{BM} (the case of the Banach-Mazur distance), we also need the following statement.

\begin{lemma}
\label{BMKZW}
If $\tau\in(0,1/4)$ and the $o$-symmetric convex bodies $K,Z\subset\R^n$ satisfy $K\subseteq Z$,  
$(1-\tau) W^n\subseteq Z$, $(1-2\tau)W^n\not\subseteq K$ and $V(Z)\leq V(W^n)$, then 
$V(K)\leq (1-\frac{\tau^n}{2^n})V(W^n)$.
\end{lemma}
\proof
Let $e_1,\ldots,e_n$ be the orthonormal basis of $\R^n$ such that the facets of $W_n$ touch $S^{n-1}$ at $\{\pm e_1,\ldots,\pm e_n\}$. Possibly reindexing $e_1,\ldots,e_n$, we may assume  for some $t>0$ that we have
\begin{eqnarray*}
t\sum_{i=1}^ne_i&\in &\partial K, \text{ and}\\
t\sum_{i=1}^n\eta_ie_i &\in &K
\mbox{ if $\eta_i\in\{-1,1\}$, $i=1,\ldots,n$, and some $\eta_i\neq 1$.}
\end{eqnarray*}
Since $(1-2\tau)W^n\not\subseteq K$, we have $t< 1-2\tau$. It follows that
\begin{eqnarray*}
({\rm int}\,K)\cap \left(\tau [0,1]^n+t\sum_{i=1}^ne_i\right)&=&\emptyset,\\
\tau [0,1]^n+t\sum_{i=1}^ne_i&\subseteq & (1-\tau) W^n\subseteq Z.
\end{eqnarray*}
Therefore
$$
V(K)\leq V(Z)-\tau^n\leq \left(1-\frac{\tau^n}{2^n}\right)V(W^n).
\mbox{ \ \proofbox}
$$

\noindent{\bf Proof  of Corollary~\ref{Zinftymustab} } We may assume that $\mu$ is not a cross measure. 
For an even isotropic measure $\mu$ and  a sufficiently small $\varepsilon>0$, we assume that 
\begin{eqnarray}
\label{Z*inftycond}
V(Z^*_{\infty}(\mu))&\ge& (1-\varepsilon) V(Z^*_{\infty}(\nu_n))\\
\label{Zinftycond}
\mbox{or \ }V(Z_{\infty}(\mu))&\le& (1+\varepsilon) V(Z_{\infty}(\nu_n)),
\end{eqnarray}
and prove that
$$
\delta_{HO}({\rm supp}\,\mu,{\rm supp}\,\nu_n)<n^{cn^3}\varepsilon^{1/3}
$$
for some absolute constant $c>0$. How small $\varepsilon$ should be is specified by 
(\ref{epsiloncond}).

According to Theorem~\ref{Zpmustab}, there exists an absolute constant $c_0>0$ such that if 
$n^{c_0n^3}\varepsilon^{1/3}<1$, then \eqref{Z*inftycond} implies that 
\begin{equation}
\label{mununassumption}
\delta_{\rm W}(\mu,\nu_n)<n^{c_0n^3}\varepsilon^{1/3},
\end{equation}
where ${\rm supp}\,\nu_n=\{\pm e_1,\ldots,\pm e_n\}$ for an orthonormal basis 
$e_1,\ldots,e_n$ of $\R^n$. In particular, $Z^*_{\infty}(\nu_n)=W^n$, and
$Z_{\infty}(\nu_n)$ is the cross polytope $C^n=[\pm e_1,\ldots,\pm e_n]$, where 
$[z_1,\ldots,z_k]$ denotes the convex hull of points $z_1,\ldots,z_k\in\R^n$. 

In the following argument, we require that 
\begin{equation}
\label{epsiloncond}
3n^26^nn! n^{c_0n^3}\varepsilon^{1/3}<\pi/4.
\end{equation}

We claim that for any $i\in\{1,\ldots,n\}$ there exists $u_i\in{\rm supp}\,\mu$ such that
\begin{equation}
\label{uiclosetoei}
\angle(u_i,e_i)\leq n^{c_0n^3}\varepsilon^{1/3}.
\end{equation}
We suppose that say for $e_1$,  we have $\angle (e_1,u)>n^{c_0n^3}\varepsilon^{1/3}$ 
for any $u\in{\rm supp}\,\mu$, and seek a contradiction. 
Naturally, also $\angle (-e_1,u)>n^{c_0n^3}\varepsilon^{1/3}$ for any $u\in{\rm supp}\,\mu$. 
We consider the function $f\in {\rm Lip}_1(S^{n-1})$ defined by
$$
f(u)=\max\left\{0,n^{c_0n^3}\varepsilon^{1/3}-\angle(u,e_1),n^{c_0n^3}\varepsilon^{1/3}-\angle(u,-e_1)
\right\}\mbox{ \ \ for $u\in S^{n-1}$}.
$$
Then we have
$$
\int_{S^{n-1}}f\,d\nu_n=n^{c_0n^3}\varepsilon^{1/3}\mbox{ \ and \ }
\int_{S^{n-1}}f\,d\mu=0,
$$
contradicting (\ref{mununassumption}), and proving (\ref{uiclosetoei}). 
Writing $\mu_0$ to denote any even measure on $S^{n-1}$ with support 
$\{\pm u_1,\ldots,\pm u_n\}$, we deduce from (\ref{uiclosetoei}) and Lemma~\ref{cubesandwich} that
\begin{equation}
\label{Z*inftyinW}
Z^*_\infty(\mu)\subseteq Z^*_\infty(\mu_0)\subseteq  e^{2n\alpha}W^n
\mbox{ \ for $\alpha=n^{c_0n^3}\varepsilon^{1/3}$.}
\end{equation}

Let $w=\sum_{i=1}^n e_i$, let $\varphi=\min\left\{\delta_{\rm H}(\mu,\nu_n),\frac{\pi}4\right\}$, and let
$u\in{\rm supp}\,\mu$ be such that $\angle(u,e_i)\geq \varphi$ and $\angle(u,-e_i)\geq \varphi$ for 
$i=1,\ldots,n$. In particular, $\varphi\in (0,\frac{\pi}4]$ as $\mu\neq \nu_n$. 
Possibly after changing the sign of some of the vectors $e_1,\ldots,e_n$, we may assume that 
$u\in {\rm pos}\,\{e_1,\ldots,e_n\}$. Let $u=(t_1,\ldots,t_n)$, where we may assume that
$$
0\leq t_1\leq \ldots \leq t_n \leq \cos \varphi.
$$

We prove that
\begin{equation}
\label{sumti}
\langle u,w\rangle \geq 1 +\frac{\varphi}3.
\end{equation}
Our task is to minimize $\langle u,w\rangle=\sum_{i=1}^nt_i$ under the conditions that
each $t_i\in[0,\cos\varphi]$ and $\sum_{i=1}^nt^2_i=1$. Solving this problem leads to
$$
\langle u,w\rangle=\sum_{i=1}^n t_i\geq \cos\varphi+\sin\varphi=\sqrt{1+\sin2\varphi}>1+\frac{\sin2\varphi}3,
$$
proving (\ref{sumti}).

First, we assume that (\ref{Z*inftycond}) holds. For the halfspace $H^+=\{x\in\R^n:\langle x, u\rangle \geq 1\}$, 
we claim that
\begin{equation}
\label{H+W}
V(H^+\cap W^n)\geq \frac{\varphi}{6^nn!}\, V(W^n).
\end{equation}
 For $i=1,\ldots,n$, let $s_i\in[0,2]$ be maximal such that $w-s_i e_i\in H^+\cap W^n$. Then we have
 $\langle w-s_i e_i, u\rangle =1$ provided $s_i<2$, thus (\ref{sumti}) yields
$$
s_i=\min\left\{2, \frac{\langle u,  w\rangle -1}{t_i}\right\}\geq \min\left\{2, \frac{\varphi}{3t_i}\right\},
$$
where we use the convention  $\frac{a}{0}=\infty$ for $a>0$. We consider two cases. 
If $\varphi=\frac{\pi}4$, then $t_i< \varphi$, and hence $s_i\geq 1/3$ for $i=1,\ldots,n$. We deduce that
$$
V(H^+\cap W^n)\geq \frac{s_1\cdots s_n}{n!}\geq \frac1{3^n n!}\geq 
\frac{\varphi}{6^nn!}\, V(W^n).
$$
If $0<\varphi<\frac{\pi}4$, then $t_n=\cos\varphi$, 
thus $t_i\leq \sin\varphi<\varphi$ for $i=1,\ldots,n-1$. In particular, 
$s_n>\frac{\varphi}3$, and $s_i>\frac13 $ for $i=1,\ldots,n-1$, and hence
$$
V(H^+\cap W^n)\geq \frac{s_1\cdots s_n}{n!}\geq \frac{\varphi}{3^n n!}=
\frac{\varphi}{6^nn!}\, V(W^n).
$$

We deduce from $2n^2\alpha<1$ (cf.~(\ref{epsiloncond})),  (\ref{Z*inftyinW}) and (\ref{H+W}) that
$$
V(Z^*_\infty(\mu))\leq e^{2n^2\alpha}V(W^n)-\frac{2\varphi}{6^nn!}\, V(W^n)
\leq \left(1+4n^2n^{c_0n^3}\varepsilon^{1/3}-\frac{2\varphi}{6^nn!}\right)V(W^n).
$$
Comparing to (\ref{Z*inftycond}) yields that
$$
\varphi<3n^2 6^nn! n^{c_0n^3}\varepsilon^{1/3},
$$
where $\delta_{\rm H}(\mu,\nu_n)=\varphi$ by (\ref{epsiloncond}).

Finally we assume (\ref{Zinftycond}). We deduce from (\ref{Z*inftyinW}) and by duality that
$$
e^{-2n\alpha}C^n\subseteq Z_\infty(\mu). 
$$
Let $T_o=[o,e^{-2n\alpha}e_1,\ldots,e^{-2n\alpha}e_n]$ and 
$T_u=[u,e^{-2n\alpha}e_1,\ldots,e^{-2n\alpha}e_n]$. Since the height of the simplex 
$T_u$ corresponding to $u$ is $n^{-1/2}(\langle u, w\rangle -e^{-2n\alpha})$, and the height of $T_o$ 
corresponding to $o$ is $n^{-1/2}e^{-2n\alpha}$,
it follows from (\ref{sumti}) that
$$
V(T_u)\geq \frac{\varphi}3\,V(T_o)= \frac{\varphi}{3\cdot 2^n}\,V(e^{-2n\alpha}C^n).
$$
Since $u\in{\rm supp}\,\mu$, we have
$$
V(Z_\infty(\mu))\geq \left(1+\frac{\varphi}{3\cdot 2^n}\right)e^{-2n^2\alpha}V(C^n).
$$
Comparing to (\ref{Zinftycond}) implies that
$$
1+\frac{\varphi}{3\cdot 2^n}\leq e^{2n^2\alpha}(1+\varepsilon)<e^{3n^2\alpha}<1+6n^2n^{c_0n^3}.
$$
We conclude $\varphi\leq 18\cdot 2^nn^2n^{c_0n^3}$, where
$\delta_{\rm H}(\mu,\nu_n)=\varphi$ by (\ref{epsiloncond}).
\proofbox

\noindent{\bf Proofs of Theorems~\ref{vol} and \ref{BM}: } Let $K$ be an origin symmetric 
convex body such that $B^n$ is the maximal volume ellipsoid contained in $K$, and suppose that 
\begin{equation}
\label{volBMcond}
\frac{S(K)^n}{V(K)^{n-1}}\geq(1-\varepsilon)\frac{S(W^n)^n}{V(W^n)^{n-1}}
\end{equation}
for a sufficiently small $\varepsilon>0$. If $C$ is a compact convex set with $B^n\subseteq C$, 
and $S_C$ is the surface area measure of $C$, then
$$
V(C)=\int_{S^{n-1}}\frac{h_C(u)}n\,d S_C(u)\geq \int_{S^{n-1}}\frac{1}n\,d S_C(u)=\frac{S(C)}n,
$$
with equality if $h_C(u)=1$ for each $u\in{\rm supp} \,S_C$.
 Therefore $V(W^n)=S(W^n)/n$ and  $V(K)\geq S(K)/n$, and hence
(\ref{volBMcond}) implies 
\begin{equation}
\label{VKVW}
V(K)\geq (1-\varepsilon)V(W^n).
\end{equation}

Let $\mu$ be a discrete even  isotropic measure satisfying 
${\rm supp}\,\mu\subseteq S^{n-1}\cap\partial K$ provided by John's Theorem. In particular,
\begin{equation}
\label{KZinfty}
\mbox{$K\subseteq Z^*_{\infty}(\mu)$ \ and \ }V(Z^*_{\infty}(\mu))\geq V(K)\geq (1-\varepsilon)V(W^n).
\end{equation}
We deduce from Corollary~\ref{Zinftymustab} that, possibly after a suitable rotation, we may assume that 
$$
\delta_{H}({\rm supp}\,\mu,{\rm supp}\,\nu_n)\leq n^{c_1n^3}\varepsilon^{\frac13}
$$
 for an absolute constant $c_1>0$. Applying now Lemma~\ref{cubesandwich}, we have
\begin{equation}
\label{cubesandwich0}
e^{- \omega\varepsilon^{\frac13}}W^n\subseteq  
Z^*_\infty(\mu)\subseteq e^{ \omega\varepsilon^{\frac13}}W^n
\end{equation}
for $\omega=n^{c_2n^3}$ and an absolute constant $c_2>0$ (assuming that $\varepsilon$ is sufficiently small).

To verify the estimate of Theorem~\ref{vol} for $\delta_{\rm vol}$, let us write 
$\delta_{\rm sym}(C,M)=V(C\Delta M)$ to denote the distance of two compact convex sets 
according to the symmetric difference metric. For example, (\ref{cubesandwich0}) yields
$$
\delta_{\rm sym}(Z^*_\infty(\mu), W^n)\leq 
\left( e^{n\omega\varepsilon^{\frac13}}-
e^{-n\omega\varepsilon^{\frac13}}\right)
2^n\leq n^{c_3n^3}\varepsilon^{\frac13}\cdot 2^n
$$
for an absolute constant $c_3>0$. We note that $V(K)\leq V(Z^*_\infty(\mu))\leq 2^n$
by K.M. Ball's Theorem~B. Hence,
$$
0\le \delta_{\rm sym}(Z^*_\infty(\mu), K)=V(Z^*_\infty(\mu))-V(K)
\leq V(Z^*_\infty(\mu))- V(W^n)+2^n\varepsilon\le 2^n\varepsilon.
$$
Let $\lambda \geq 1$ be such that $V(\lambda K)=2^n$,
and hence $V(\lambda K)-V(K)\leq \varepsilon\cdot 2^n$ according to (\ref{KZinfty}).
We conlude that
\begin{eqnarray*}
\delta_{\rm vol}(K,W^n)&\leq &2^{-n}\delta_{\rm sym}(\lambda K,W^n)\\
&\leq&
2^{-n}(\delta_{\rm sym}(\lambda K,K)+
\delta_{\rm sym}(K, Z^*_\infty(\mu))+
\delta_{\rm sym}(Z^*_\infty(\mu), W^n))\\
&\leq & n^{c_4n^3}\varepsilon^{\frac13},
\end{eqnarray*}
for an absolute constant $c_4>0$, and this completes the proof of Theorem~\ref{vol}.

Let us turn to the estimate of Theorem~\ref{BM} for $\delta_{\rm BM}$. Let
$\delta_{\rm BM}(K,W^n)\geq \alpha$ for some $\alpha\in(0,1)$. If 
\begin{equation}
\label{BMZcond}
e^{-\frac{\alpha}5}W^n \subseteq Z^*_\infty(\mu)\subseteq e^{\frac{\alpha}5}W^n,
\end{equation}
then $\delta_{\rm BM}(K,W^n)\geq \alpha$ implies that
$e^{-\frac{4\alpha}5}W^n \not\subseteq K$, and hence
$(1-\frac{2\alpha}5)W^n\not\subseteq K$. On the other hand,
$(1-\frac{\alpha}5)W^n\subseteq Z^*_\infty(\mu)$, thus Lemma~\ref{BMKZW} yields
\begin{equation}
\label{BMVK1}
V(K)\leq \left(1-\frac{\alpha^n}{10^n}\right)V(W^n).
\end{equation}

Finally, we assume that (\ref{BMZcond}) does not hold. Since (\ref{VKVW}) leads to
(\ref{cubesandwich0}), we have $V(K)< (1-\varepsilon)V(W^n)$ provided
 $\frac{\alpha}5=\omega\varepsilon^{\frac13}$. In other words,
\begin{equation}
\label{BMVK2}
V(K)\leq \left(1-\frac{\alpha^3}{125\omega^3}\right)V(W^n)
\end{equation}
where $\frac1{125\omega^3}\geq n^{-c_5n^3}$ for an absolute constant $c_5>0$.
Combining 
(\ref{BMVK1}) and (\ref{BMVK2}) proves Theorem~\ref{BM}.
\proofbox

\section{Proof of Theorem~\ref{volBMplanar}}

In this section, we prove Theorem~\ref{volBMplanar}, which is the $2$-dimensional 
(sharper) version of Theorems~\ref{vol} and \ref{BM}. The idea of our proof is 
essentially the one given by F.~Behrend \cite{Behrend1937}. 
As before, let $[x_1,\ldots,x_k]$ denote the convex hull of the points 
$x_1,\ldots,x_k\in\R^2$. For the origin symmetric convex body $K\subseteq\R^2$ and 
$u\in \R^2\setminus\{o\}$, we write $H(K,u)$ to denote the supporting line with 
exterior normal $u$, and $H(K,u)^-$ to denote the corresponding halfplane containing $K$.

Let $\varepsilon\in[0,\frac12)$. Let $K$ be a planar origin symmetric convex body 
which has a square as an inscribed 
parallelogram of maximum area. Suppose that  
\begin{equation}
\label{planarcond}
\frac{S(K)^2}{V(K)}\geq (1-\varepsilon)\frac{S(W^2)^2}{V(W^2)}.
\end{equation}
Then we  prove that
\begin{eqnarray}
\label{planarvol}
\delta_{\rm vol}(K,W^2)&\leq & 54\varepsilon \text{ and} \\
\label{planarBM}
\delta_{\rm BM}(K,W^2)&\leq & 18\varepsilon. 
\end{eqnarray}
Let $u_1,u_2$ denote the standard basis of $\R^2$. We may assume that $W^2=[-1,1]^2$ is 
a  parallelogram of largest area contained in $K$, and hence $p_i\in\partial K\cap H(K,p_i)$ 
holds for the vertices $p_1=u_2+u_1$ and 
$p_2=u_2-u_1$ of $W^2$. It also follows that
\begin{equation}
\label{Kinsquare}
K\subseteq\bigcap_{i=1}^2H(K,\pm p_i)^-=[\pm 2u_1, \pm 2u_2].
\end{equation}
Let $q_i\in \partial K\cap H(K,u_i)$ for $i=1,2$. In particular, (\ref{Kinsquare}) yields
\begin{eqnarray*}
q_1&= & (1+t_1,s_1) \mbox{ where $t_1\in[0,1]$ and $|s_1|\leq 1-t_1$,}\\
q_2&= & (s_2,1+t_2) \mbox{ where $t_2\in[0,1]$ and $|s_2|\leq 1-t_2$.}
\end{eqnarray*}
Since $K$ contains the parallelogram $P=[\pm q_1,\pm q_2]$, we have
\begin{eqnarray*}
V(W^2)&\geq & V(P)=2|\det[q_1,q_2]|=2[(1+t_1)(1+t_2)-s_1s_2]\\
&\geq&2[(1+t_1)(1+t_2)-(1-t_1)(1-t_2)]=4(t_1+t_2),
\end{eqnarray*}
and hence
$$
t=\frac{t_1+t_2}2\leq \frac12.
$$
We approximate $K$ by suitable polygons to obtain
\begin{equation}
\label{QKM}
W^2\subseteq Q\subseteq K\subseteq M\subseteq (1+t)W^2,
\end{equation}
where
\begin{eqnarray*}
M&=&\left(\bigcap_{i=1}^2H(K,\pm u_i)^-\right)\bigcap\left(\bigcap_{i=1}^2H(K,\pm p_i)^-\right)
\mbox{ with $S(M)=(1+(\sqrt{2}-1)t)S(W^2)$},\\
Q&=&[\pm p_1,\pm p_2,\pm q_1,\pm q_2]\mbox{ with $V(Q)=(1+t)V(W^2)$}.
\end{eqnarray*}
We deduce from (\ref{planarcond}) and (\ref{QKM}) that
$$
(1-\varepsilon)\frac{S(W^2)^2}{V(W^2)}\leq
\frac{S(K)^2}{V(K)}\leq\frac{S(M)^2}{V(Q)}
=\frac{(1+(\sqrt{2}-1)t)^2S(W^2)^2}{(1+t)V(W^2)}.
$$
Since $\frac{1-t}{1+t}\geq\frac13$ by $t\leq \frac12$, we have
\begin{equation}
\label{epsilont}
\varepsilon\geq 1-\frac{(1+(\sqrt{2}-1)t)^2}{1+t}=\frac{(3-2\sqrt{2})t(1-t)}{1+t}\geq
\frac{(3-2\sqrt{2})t}{3}\geq\frac{t}{18}.
\end{equation}
Therefore combining (\ref{QKM}) and (\ref{epsilont}) leads to
$$
\delta_{\rm BM}(K,W^2)\leq \log(1+t)\leq t\leq 18\varepsilon,
$$
and combining (\ref{QKM}) and (\ref{epsilont}) with an elementary argument leads to
$$
\delta_{\rm vol}(K,W^2)\leq (1+t)^2-1\leq 3t\leq 54\varepsilon.
$$
We conclude (\ref{planarvol}) and (\ref{planarBM}), and in turn Theorem~\ref{volBMplanar}.


\section{Even isotropic measures on $S^1$}

The goal of this section is to prove the following improvement of Corollary~\ref{Zinftymustab} if $n=2$.

\begin{theo}
\label{Zinftymustab2}
If  $\mu$ is an even isotropic measure on $S^1$, $\varepsilon\in(0,1)$, and 
$\delta_{HO}({\rm supp}\,\mu,{\rm supp}\,\nu_2)\geq  \varepsilon$, then
\begin{eqnarray*}
V(Z_{\infty}(\mu))&\ge& (1+ 0.25 \,\varepsilon) V(Z_{\infty}(\nu_2)),\\
V(Z^*_{\infty}(\mu))&\le& (1- 0.1\, \varepsilon) V(Z^*_{\infty}(\nu_2)).
\end{eqnarray*}
\end{theo}

We call a compact, symmetric set $X\subseteq S^1$ proper if for each $v\in S^1$ there is some 
$u\in X$ such that $\angle(u,v)\le \pi/4$. A compact, symmetric set $X\subseteq S^1$ is proper 
if and only if the angle between consecutive points of $X$ on $S^1$ is at most $\pi/2$. For a closed 
set $X\subseteq S^1$ we define 
$$
d_0(X)=\min\{\delta_H(X,\rho\{\pm e_1,\pm e_2\}):\rho\in\text{SO}(2)\},
$$
where $e_1,e_2$ is an orthonormal basis of $\R^2$. If $X$ is proper, then $d_0(X)\le \pi/4$. 

Note that if $\mu$ is an even isotropic measure on $S^1$, then Claim~\ref{isotropiccap} shows that 
the support of $\mu$ is a proper set. 

\begin{lemma}\label{existlem}
If $X\subseteq S^1$ is proper, $\eta\in(0,\pi/4)$ and $d_0(X)\ge \eta$, then there are $u,v\in X$ such that 
$\eta\le \angle (u,v)\le \frac{\pi}{2}-\eta$.
\end{lemma}

\proof Assume that for any pair $u,v\in X$ either $\angle (u,v)<\eta$ or  $\angle (u,v)>\frac{\pi}{2}-\eta$. 
Let $u_1\in X$ be arbitrary. Then there is no $v\in X$ such that $\angle (u,v)\in [\eta,\frac{\pi}{2}-\eta]$. 
The same is true for $-u_1\in X$. Let $\bar u_1\in S^1\cap u_1^\perp$. Then there is some $u_2\in X$ with 
$\angle(\bar u_1,u_2)<\eta$. Since $X$ is closed and symmetric, we conclude that $d_0(X)<\eta$, a contradiction. 
\proofbox

We turn to the proof of Theorem \ref{Zinftymustab2} and start with the second assertion. Let the assumptions be fulfilled. 
By an approximation argument (see Barthe \cite{Bar04}), we can assume that $\mu$ is discrete. In the following, we 
use property (P) which states that for $0\le \beta\le \alpha<\pi/2$ the function 
$$F(t):=\tan\left(\frac{\alpha+t}{2}\right)+\tan\left(\frac{\beta-t}{2}\right), 
\qquad t\in[0,\min\{\beta,\tfrac{\pi}{2}-\alpha\}],
$$ 
is strictly increasing. Applying (P) repeatedly to angles between consecutive 
vectors of $\text{supp}\, \,\mu$, Lemma \ref{existlem} and symmetry, we obtain
$$
V(Z^*_{\infty}(\mu))\le 2\left(\tan\left(\frac{\alpha}{2}\right)+\tan\left(\frac{\pi}{4}-\frac{\alpha}{2}\right)+
\tan\left(\frac{\pi}{4}\right)\right)
$$
for some $\alpha\in [\varepsilon,\frac{\pi}{2}-\varepsilon]$. Since 
$$
\tan\left(\frac{\alpha}{2}\right)+\tan\left(\frac{\pi}{4}-\frac{\alpha}{2}\right)=2(1+\sin\alpha+\cos\alpha)^{-1}
$$
and 
\begin{equation}\label{anlem}
\sin\alpha+\cos\alpha\ge 1+0.5\, \varepsilon \qquad \text{for } \alpha\in[\varepsilon,\tfrac{\pi}{2}-\varepsilon]
\end{equation}
with $\varepsilon\in (0,\pi/4)$, we obtain
$$
V(Z^*_{\infty}(\mu))\le 2\left(\frac{1}{1+0.25\, \varepsilon}+1\right)<4\left(1-0.1\, \varepsilon\right),
$$
which proves the second assertion.

For the first assertion, we argue similarly. Here we use the fact that for $0\le \beta\le \alpha<\pi/2$ the function $G(t)=
\sin(\alpha+t)+\sin(\beta-t)$, $t\in[0,\min\{\beta,\tfrac{\pi}{2}-\alpha\}]$, is strictly decreasing. Thus we obtain 
$$
V(Z_{\infty}(\mu))\ge \sin(\alpha)+\sin\left(\frac{\pi}{2}-\alpha\right)+\sin\left(\frac{\pi}{2}\right)=\sin\alpha+\cos\alpha+1 
$$
for some $\alpha\in [\varepsilon,\frac{\pi}{2}-\varepsilon]$. Now the first assertion follows from \eqref{anlem}.
\proofbox

 
\mbox{ }

\section*{Acknowledgements}
K.J. B\"or\"oczky and F. Fodor are supported by 
National Research, Development and Innovation Office -- NKFIH grant 116451, and K.J. B\"or\"oczky is also supported by grant 109789. 

F. Fodor wishes to thank the Alfr\'ed R\'enyi Institute of Mathematics of the Hungarian Academy of Sciences where part of his work was done while he was a visiting researcher. 

D. Hug is supported by DFG grants FOR 1548 and HU 1874/4-2.

\bigskip

Authors' addresses:

\medskip

K\'aroly J. B\"or\"oczky, MTA Alfr\'ed R\'enyi Institute of Mathematics, Hungarian Academy of Sciences, Re\'altanoda u. 13-15, 1053 Budapest, Hungary. E-mail: carlos@renyi.hu

Ferenc Fodor, Department of Geometry, Bolyai Institute, University of Szeged, Aradi v\'ertan\'uk tere 1, 6720 Szeged, Hungary. E-mail: fodorf\@math.u-szeged.hu

Daniel Hug, Karlsruhe Institute of Technology (KIT), D-76128 Karlsruhe, Germany. E-mail: daniel.hug@kit.edu

\end{document}